\documentclass[12pt]{article}
\usepackage{amsmath,amsthm,amssymb}
\usepackage{amsfonts}
\usepackage{latexsym} 
\usepackage{graphicx}
\usepackage{proba}

\newtheorem{lem}{Lemma}
\newtheorem{pro}{Proposition}
\newtheorem{cor}{Corollary}
\newtheorem{remark}{Remark}

\newcommand{\QQ}{{\kern.24em\vrule width.04em height1.4ex%
                 depth-.05ex\kern-.26em\mathsf Q}}
\newcommand{\CC}{{\kern.24em\vrule width.04em height1.4ex%
                 depth-.05ex\kern-.26em\mathsf C}}

\def\SC#1{{\mathcal{#1}}}

\def\BB#1{{\mathbb{#1}}}

\title{\Large\bfseries
Shrinkage estimation of a mean matrix of a multivariate complex normal distribution}
\author{
Yoshihiko Konno
\thanks{
\hskip -18pt 2000 {\it Mathematics Subject Classification.} Primary: 62H12, Secondary:  62F10.
\hfill\break
\hskip 2cm{\it Key words and phrases.} unbiased estimation of risk, integration-by-parts identities, complex Wishart distribution 
}
}
\begin{document}
\maketitle
\begin{abstract}\noindent
The problem of estimating a mean matrix of a multivariate complex normal distribution with an unknown covariance matrix is considered under an invariant loss function. By using complex versions of the Stein identity, the Stein-Haff identity, and calculus on eigenvalues, a formula is obtained for an unbiased estimate of the risk of an invariant class of estimators, from which several minimax shrinkage estimators are constructed. 
\end{abstract}
\section{Introduction}
The multivariate complex normal and complex Wishart distributions were first explored in Goodman~\cite{Goodman1963}, and followed by Khatri~\cite{Khatri1965}. These models play an important role in signal processing methods. See Kay~\cite{Kay1993} for the need of complex data models  and DoGond\v{z}i\'{c} and Neborai~\cite{DN2003} for a unified approach based on complex GMANOVA models to analyze and extend signal processing models. See 
Ratnarajah  {\it et al.}~\cite{Ratnarajah2005}, Micheas {\it et al.}~\cite{Micheas2006}, and D\'{i}az-Garc\'{i}a and Gutierr\'{e}z-J\'{a}imez~\cite{Diaz-Garcia2011} for recent development of complex data model.  
Lillest{\o}l~\cite{Lillestol1977} first investigated Stein-like shrinkage methods on simultaneous estimation of a mean vector of the complex normal model. However, shrinkage methods for these models have received less attention so far, although  it is important to develop these methods beyond the maximum likelihood estimator of estimating the unknown signals in the multivariate complex normal distribution. The goal of this paper is to show how certain decision theoretical results concerning the problem of estimating a mean matrix of the real normal distribution can be extended to the complex multivariate normal case.   
\par
In this paper, we consider the problem of estimating an $ m \times p $ unknown constant complex matrix $ \mathbf{\Xi}$ that is observed with additive complex normal random errors in a decision theoretic set-up.  
Our observations are an $ m \times p $ data matrix $ \mathbf{Z} $ and a 
$ p \times p $ positive definite Hermitian matrix $ \mathbf{S}$, which is represented as 
\begin{eqnarray}
  \begin{array}{l}
  \mathbf{Z}:\, m \times p 
  \sim
  \BB{C} N _{ m \times p } ( 
    \mathbf{\Xi}
    ,\,
    \mathbf{K}  \otimes \mathbf{\Sigma} 
  ),
  \\
  \mathbf{S}:\, p \times p 
  \sim
  \BB{C} W _p (\mathbf{\Sigma},\,n )
  \qquad
  \mbox{with $\mathbf{Z}$ and $\mathbf{S}$ independent},
  \end{array}
  \label{eq:0}
\end{eqnarray} 
where $ n > p$, $ \mathbf{\Sigma}$ is a  $ p \times p $ positive definite Hermitian constant matrix, and $ \mathbf{K} $ is an $ m \times m $  positive definite Hermitian constant matrix. 
Here we assume that $ \mathbf{\Xi} $ and $ \mathbf{\Sigma} $ are unknown although we assume that $ \mathbf{K} $ is known.  
Furthermore $ \BB{C} N _{ m \times p } ( 
    \mathbf{\Xi}
    ,\,
    \mathbf{K}  \otimes \mathbf{\Sigma} 
  ) $ and 
  $   \BB{C} W _p (\mathbf{\Sigma},\,n ) $ stand for a matrix-variate complex normal distribution with the mean matrix $ \mathbf{\Xi} $ and the covariance matrix 
  $ \mathbf{K}  \otimes \mathbf{\Sigma} $ and a complex Wishart distribution with the degree of freedom $n $ and the parameters $  \mathbf{\Sigma} $, respectively. 
In other words,  
the model $(\ref{eq:0})$ means that the density of $ \mathbf{K} ^{-1/2} \mathbf{Z}=:\widetilde{\mathbf{Z}} $ with respect to the Lebesgue measure on 
$ \BB{C} ^{ m \times p } $ is given as 
\begin{equation}
  \pi ^{ - mp }
  {\rm Det} ( \mathbf{\Sigma} ) ^{ - m }
  \exp \bigl \{
    - {\rm Tr} \bigl (
      ( \mathbf{z} - \widetilde{\mathbf{\Xi}} ) \mathbf{\Sigma}^{-1 } 
      ( \mathbf{z} - \widetilde{\mathbf{\Xi}} ) ^*
    \bigr )
  \bigr \},
  \qquad
  \mathbf{z} \in \BB{C} ^{ m \times p }, 
  \nonumber
\end{equation}
where $ \widetilde{\mathbf{\Xi}} = \mathbf{K} ^{-1/2} \mathbf{\Xi}$, while 
the density of $ \mathbf{S} $ with respect to the Lebesgue measure on 
$\BB{C} ^{ p \times p } _+ $ is given by 
\begin{equation}
  \frac{
    {\rm Det\,}( \mathbf{s} ) ^{ n-p} \exp  
    \bigl (
      - {\rm Tr\,} ( \mathbf{s} \mathbf{\Sigma} ^{-1} )
    \bigr )
  }{
    {\rm Det\,}(\mathbf{\Sigma})^n 
    \pi ^{ p(p-1)/2} \Pi _{ k = 1 } ^p
    \Gamma ( n + 1 -k )
  },
  \qquad
  \mathbf{s} \in \BB{C} _+ ^{ p \times p }. 
  \label{eq:den2}
\end{equation}
Here $ \Gamma (\,\cdot\,) $ is the usual Gamma function, $ {\rm Tr\,}(\,\cdot\,) $ and $ {\rm Det\,}(\,\cdot\,) $ denote the trace and determinant of a square matrix, and the superscript  
''${}^*$'' means the complex conjugate transpose of a matrix. Furthermore 
$ \BB{C} ^{ m \times p } $ and $ \BB{C} ^{ p \times p } _+ $ stand for the sets of all $ m \times p $ complex matrices and of all $ p \times p $ positive definite Hermitian complex matrices, respectively.  
\par
Based on $ ( \mathbf{Z},\,\mathbf{S} ) $ we consider the problem of estimating 
the mean matrix $ \mathbf{\Xi} $ with respect to a loss function 
\[
  \SC{L} ( \widehat{\mathbf{\Xi}},\,(\mathbf{\Xi},\,\mathbf{\Sigma} ) )
  =
  {\rm Tr\,} \{
    \mathbf{\Sigma} ^{ - 1 }
    (
      \widehat{\mathbf{\Xi}} - \mathbf{\Xi} 
    ) ^*
    \mathbf{K} ^{-1}
    (
      \widehat{\mathbf{\Xi}} - \mathbf{\Xi} 
    )
  \},
\]
where an $ m \times p $ random matrix $ \widehat{\mathbf{\Xi}} $ is an estimator of $ \mathbf{\Xi} $. 
The risk function corresponding to this loss function is 
\[
  \SC{R} ( \widehat{\mathbf{\Xi}},\,(\mathbf{\Xi},\,\mathbf{\Sigma} ) )
  =
  \BB{E} [ 
    \SC{L} ( \widehat{\mathbf{\Xi}},\,(\mathbf{\Xi},\,\mathbf{\Sigma} ) ) 
  ],
\]
where the expectation above is taken with respect to the joint distribution of 
$ (\mathbf{Z},\,\mathbf{S})$. 
\par
This estimation problem is important since it is a prototype of estimating the regression matrix of a complex MANOVA model and of predicting multivariate responses in a linear regression complex model. 
We extend a large body of the results obtained by 
Efron and Morris~\cite{EF1976}, Bilodeau and Kariya~\cite{BK1989}, Kariya {\it et al.}~\cite{KKS1999}, Konno~\cite{Konno1991}, and van der Merwe and Zidek~\cite{MZ1980} in the multivariate real normal set-up to the complex normal set-up $(\ref{eq:0})$. 
The results  in the real normal model were obtained by extensive use of the integration by parts approach, known as the Stein identity derived by Stein~\cite{Stein1973, Stein1981}, and the Stein-Haff identity by Stein~\cite{Stein1977} and Haff~\cite{Haff1979a,Haff1979b}. In addition to these identities, 
the eigenvalue calculus, developed by Loh~\cite{Loh1988,Loh1991a,Loh1991b}, Konno~\cite{Konno1991}, and Kariya {\it et al.}~\cite{KKS1999}, is important to the development for a systematic search for shrinkage estimators. 
We extend these approaches to the complex normal set-up. 
The Stein identity for the multivariate complex normal is easily derived by 
using an isomorphism between real and complex variables stated in Andersen {\it et al.}~\cite{Anderson1985} while the Stein-Haff identity was extended to the 
complex Wishart distribution by Svensson and Lundberg~\cite{Svensson2004}. 
These identities and the eigenvalue calculus for the complex matrix developed in this paper are exploited to establish a systematic search for shrinkage estimators for the model $(\ref{eq:0})$, which includes the FICYREG  estimator of van der Merwe and Zidek~\cite{MZ1980}. 
\par
Shrinkage methods for estimating the regression matrix in a  multivariate linear regression model have been extensively investigated to overcome the shortcomings of the ordinary least squares estimator. The literature includes 
Brown and Zidek~\cite{BZ1980,BZ1982}, Dempster~\cite{Dempster1977}, and 
van der Merwe and Zidek~\cite{MZ1980}. Later
 Breiman and Friedman~\cite{BF1997} proposed to predict a future observation by a ridge-type shrinkage estimator in order to use information of correlated variables. See also 
Bilodeau~\cite{Bilodeau2000}, Oman~\cite{Oman2002}, and Srivastava and Solanky~\cite{Srivastava2003} for further investigation on this problem. 
As mentioned in Srivastava and Solanky~\cite{Srivastava2003}, 
we can use minimax estimators to construct better predictors in order to overcome shortcomings of the predictor based on the least squares estimator. This shows that  the results obtained in this paper can be immediately applied to the problem of predicting a future observation in a  multivariate linear model for complex data. 
\par
The remaining parts of this papers are organized as follows. 
In Section~\ref{preliminaries}, we state some notation and the integration by parts formulae. In Section~\ref{known}, we develop shrinkage estimators for the known covariance case, which is an  extension to the results obtained in Stein~\cite{Stein1977} and Zheng~\cite{Zheng1986a,Zheng1986b}. In Section~\ref{unknown}, we obtain unbiased risk estimate for invariant estimators, from which several shrinkage estimators are derived. In the Appendix, the results on  eigenvalue calculus for the complex set-up and their proofs are developed.

\section{Preliminaries: Notation and Basic identities}\label{preliminaries}
This section first presents some notation used throughout this paper. 
Next we introduce integration by parts formulae, complex versions of 
the Stein identity and the Stein-Haff identity, which play vital roles in obtaining unbiased risk estimate in Sections~\ref{known} and \ref{unknown}. 
\subsection{Notation}\label{notation}
Let $\BB{R}$ and $\BB{C}$ denote the field of real and complex numbers, respectively. We represent any element $ c \in \BB{C} $ as 
$ c = a + \sqrt{-1}\, b $, where $ a,\,b \in \BB{R}$. We also denote the real and imaginary parts of $ c $ by ${\rm Re\,} c $ and $ {\rm Im\,} c $, 
respectively. In particular we denote by $ \BB{R} _+ $ the set of all positive real numbers. The conjugate of a complex number $ c $ is given by 
$ \bar c := a - \sqrt{ -1 }\, b $. We define by $ \BB{R} ^p $ and $ \BB{C} ^p$ the sets of all $p$-tuples of real and complex numbers, respectively. 
We set $ \BB{R} _{ > } ^p = \{ ( \ell _1,\,\ell _2,\,\ldots,\,\ell _p ) \in \BB{R} ^p:\, \ell _1 > \ell _2 > \cdots > \ell _p > 0 \} $. 
In this paper, these tuples are represented as columns. The sets of all $ m \times p $  matrices of real and complex entries are 
denoted by $ \BB{R} ^{ m \times p } $ and $ \BB{C} ^{ m \times p } $, respectively. 
The transpose and the conjugate of $ \mathbf{C} $ are denoted by $ \mathbf{C} ^\prime$ and  $ \overline{\mathbf{C}}$, respectively. 
Furthermore the conjugate transpose of an $ m \times p $ matrix 
$ \mathbf{C} \in \BB{C} ^{ m \times p } $ are denoted by $ \mathbf{C} ^* = \overline{\mathbf{C}^\prime }$. 
The set of $ p \times p $ Hermitian positive definite matrices is denoted by 
$ \BB{C} ^{ p \times p } _+$. For any $ \mathbf{c} = \mathbf{a} + \sqrt{-1}\, \mathbf{b} \in \BB{C} ^p\,( \mathbf{a},\,\mathbf{b} \in \BB{R} ^p)$, we denote by 
$ [ \mathbf{c} ] $  a $2p$-dimensional real vector $(\mathbf{a} ^\prime,\,\mathbf{b} ^\prime) ^\prime$. 
For a positive integer q and real numbers $ a _1,\,a _2,\,\ldots,\,a _q $, 
$ {\rm Diag} ( a _1,\,a _2,\,\ldots,\,a _q ) $ denotes a $ q \times q $ diagonal matrix with the $i$-th diagonal element $ a _i\,(i = 1,\,2,\,\ldots,\,q )$. 
For an $m \times p$ complex matrix 
$ \mathbf{C} = \mathbf{A} + \sqrt{-1}\, \mathbf{B}\,(\mathbf{A},\,\mathbf{B} \in \BB{R} ^{ m \times p})$, we denote by 
$ \{ \mathbf{C} \} $ a $2m \times 2p$ real matrix  
\[
  \left (
    \begin{array}{cc}
      \mathbf{A} & - \mathbf{B}
      \\
      \mathbf{B} & \mathbf{A}
    \end{array}
  \right ). 
\]
\par
Let $ g (x,\,y ) $ be a real-valued function on an open set $ U \in \BB{R} ^2 $. We say that $ g $ is differentiable if 
$ \partial g / \partial x$ and $ \partial g / \partial y $ exist on 
$ U$. 
Let $ u,\,v $ be  real-valued functions on an open set $ U \in \BB{R} ^2 $. 
A function $ g := u + \sqrt{-1}\,v $ is called differentiable if 
$ u,\,v $ are differentiable. 
For $ z = x + \sqrt{ - 1 }\, y\, (x,\,y \in \BB{R}) $ and differentiable function  $ g ( z ) = u ( z ) + \sqrt{ - 1 }\, v ( z ) $, 
we define  
\begin{eqnarray*}
  \frac{ \partial  }{ \partial z } g
  &=&
  \frac{ 1 }{ 2 } 
  \left (
    \frac{ \partial }{ \partial x } 
    - \sqrt{-1}\, 
    \frac{ \partial }{ \partial y }
  \right ) g
  =
  \frac{ 1 }{ 2 } 
  \left (
    \frac{ \partial u }{ \partial x } 
    +
    \frac{ \partial v }{ \partial y }
  \right )
  +
  \frac{ \sqrt{-1} }{ 2 } 
  \left (
    \frac{ \partial v}{ \partial x } 
    -  
    \frac{ \partial u }{ \partial y }
  \right ),
  \\
  \frac{ \partial }{ \partial \bar z } g
  &=&
  \frac{ 1 }{ 2 } 
  \left (
    \frac{ \partial }{ \partial x } 
    + \sqrt{-1}\, 
    \frac{ \partial }{ \partial y }
  \right ) g
  =
  \frac{ 1 }{ 2 } 
  \left (
    \frac{ \partial u }{ \partial x } 
    -
    \frac{ \partial v }{ \partial y }
  \right )
  +
  \frac{ \sqrt{-1} }{ 2 } 
  \left (
    \frac{ \partial v}{ \partial x } 
    +  
    \frac{ \partial u }{ \partial y }
  \right ).  
\end{eqnarray*}
It is  checked directly that 
\begin{eqnarray*}
  \frac{ \partial }{ \partial z} z = 1,
  \qquad
   \frac{ \partial }{ \partial z} \bar z = 0,
   \qquad
    \frac{ \partial }{ \partial \bar z } z = 0,
    \qquad
     \frac{ \partial }{ \partial \bar z} \bar z = 1.
\end{eqnarray*}
If $ g $ is differentiable, then 
\begin{equation}
  \overline{
    \frac{
      \partial
    }{
      \partial z
    } 
    g
  }
  =
  \frac{ 
    \partial
  }{
    \partial \bar z
  }
  \bar g.
  \label{eq:00}
\end{equation}
\par
Let $ \mathbf{G} = ( g _{ ij } ) _{ i = 1,\,2,\,\ldots,\,m,\,j = 1,\,2,\,\ldots,\,p} $ be an $ m \times p $ matrix, where 
$ g _{ij} $'s are complex-valued differentiable functions on $ \BB{C} ^{ m \times p} $. 
For $ \mathbf{z} = ( z _{ij} )  _{i = 1,\,2,\,\ldots,\,m,\,j = 1,\,2,\,\ldots,\,p}  \in \BB{C} ^{ m \times p }$, we set 
\[
  \nabla _{\mathbf{z}} 
  = 
  \left (
    \frac{ \partial }{ \partial z _{ij} }
  \right )  _{i = 1,\,2,\,\ldots,\,m,\,j = 1,\,2,\,\ldots,\,p},  
\]
and we define
\[
  {\rm Re} ({\rm Tr\,} (\nabla _z ^\prime \mathbf{G})) 
  =
  {\rm Tr\,} ( {\rm Re}  (\nabla _z ^\prime \mathbf{G})) 
  =
  \sum _{ j = 1 } ^p \sum _{ i = 1 } ^m 
  \left \{
  \frac{ \partial ({\rm Re\,} g _{ ij } )}{ \partial ({\rm Re\,} z _{ij} )}
  +
  \frac{ \partial ({\rm Im\,} g _{ ij } )}{ \partial ({\rm Im\,} z _{ij} )}
  \right \}
  . 
\]
\subsection{Complex normal distributions and the Stein identity}
Recall that a  $ p \times 1 $ complex random vector $Z$ is said to have a $p$-variate complex normal distribution with a mean vector $ \mathbf{\theta} \in \BB{C} ^p$ and a covariance matrix $\mathbf{\Sigma} \in \BB{C} _+ ^{ p \times p } $ if the density of $ Z $ with respect to  Lebesgue measure on $\BB{C} ^p$ is given as 
\[
  f _{ Z } ( \mathbf{z} )
  =
  \frac{ 1 } { \pi ^p}
  {\rm Det} (\mathbf{\Sigma} ) ^{ -1 }
  \exp 
  \{
    - ( \mathbf{z} - \mathbf{\theta} ) ^*
      \mathbf{\Sigma} ^{ -1 }
      ( \mathbf{z} - \mathbf{\theta} )
  \},
  \qquad
  \mathbf{z} \in \BB{C} ^p. 
\]
We use the notation 
$ Z\,\sim\,\BB{C} N _p ( \mathbf{\theta},\,\mathbf{\Sigma} ) $ for this. 
\begin{lem}\label{stein}
Let $ Z $ be a  $ p \times 1 $ complex random vector having 
$ \BB{C} N _p ( \mathbf{\theta},\,\mathbf{\Sigma} )$ and let 
$ \mathbf{g} = ( g _1,\,g _2,\,\ldots,\,g _p ) :\, \BB{C}^p \to \BB{C} ^p $ 
be differentiable with 
\[
  \BB{E} \biggl |
    \frac{ 
      \partial\, ( {\rm Re\,} g _i ) 
    }{ 
      \partial\, ( {\rm Re\,} z _i ) 
    }
  \biggr | _{ \mathbf{z} = Z }
  < \infty,
  \qquad
   \BB{E} \biggl |
    \frac{ 
      \partial\, ( {\rm Im\,} g _i )  
    }{ 
      \partial\, ( {\rm Im\,} z _i ) 
    }
  \biggr | _{ \mathbf{z} = Z }
  < \infty,
  \qquad
  i = 1,\,2,\,\ldots,\,p. 
\]
Then we have 
\begin{eqnarray*}
  \BB{E}  [
    ( Z - \mathbf{\theta} ) ^*
    \mathbf{\Sigma} ^{-1}
    \mathbf{g} ( Z)
    +
    \mathbf{g} ^* ( Z)
     \mathbf{\Sigma} ^{-1}
    ( Z - \mathbf{\theta} )
  ]
  &=&
  \BB{E} \left [
    \sum _{ i = 1 } ^p
    \left \{
     \frac{ 
      \partial\, ( {\rm Re\,} g _i ) 
    }{ 
      \partial\, ( {\rm Re\,} z _i ) 
    }
    +
     \frac{ 
      \partial\, ( {\rm Im\,} g _i )  
    }{ 
      \partial\, ( {\rm Im\,} z _i ) 
    }
    \right \}
    \bigg |_{ \mathbf{z} = Z}
  \right ].
\end{eqnarray*}
\end{lem}
\par\noindent
{\bf Proof}. 
Note that 
\[
   ( Z - \mathbf{\theta} ) ^*
    \mathbf{\Sigma} ^{-1}
    \mathbf{g} ( Z)
    +
    \mathbf{g} ^* ( Z)
     \mathbf{\Sigma} ^{-1}
    ( Z - \mathbf{\theta} )
    =
  2 [ Z - \mathbf{\theta} ] ^\prime
  \{ \mathbf{\Sigma} ^{-1} \}
  [ \mathbf{g} ]  
\]
and that 
$  Z\,\sim\,\BB{C} N _p ( \mathbf{\theta},\,\mathbf{\Sigma} ) $ 
if and only if 
$  [ Z ]\,\sim\, N _{ 2p } ( [\mathbf{\theta}],\,(1/2)\{\mathbf{\Sigma} \} )$, a $2p$-variate multivariate real normal distribution with a $ 2p \times 1 $ mean vector $  [\mathbf{\theta}] $ and a $ 2p \times 2p$ positive definite covariance matrix $ (1/2)\{\mathbf{\Sigma} \}$. 
By the Stein identity on a multivariate real normal distribution[see Stein~\cite{Stein1981}], we have 
\begin{eqnarray*}
  \BB{E} \bigl \{ [ Z - \mathbf{\theta} ] ^\prime
  \{ (1/2) \mathbf{\Sigma} \}^{-1} 
  [ \mathbf{g} ] \bigr \}
  =  
   \BB{E} \left \{
    \sum _{ i = 1 } ^p
    \left \{
     \frac{ 
      \partial\, ( {\rm Re\,} g _i ) 
    }{ 
      \partial\, ( {\rm Re\,} z _i ) 
    }
    +
     \frac{ 
      \partial\, ( {\rm Im\,} g _i )  
    }{ 
      \partial\, ( {\rm Im\,} z _i ) 
    }
    \right \}
    \bigg |_{ \mathbf{z} = Z}
  \right \},
\end{eqnarray*}
which completes the proof. 
\hfill $\Box$

\subsection{Complex Wishart distributions and the Stein-Haff identity}
Assume that a $ p \times p $ Hermitian positive definite matrix 
$ \mathbf{S} $  has a complex Wishart distribution 
$ \BB{C} W _p(\mathbf{\Sigma},\,n )$  with the density function $(\ref{eq:den2})$. 
Let 
$ \mathbf{G} ( \mathbf{S} ) $ be a $ p \times p $ matrix,  
the  
$ (i,\,j ) $ element $ g _{ ij } ( \mathbf{S} ) $ of which is a complex-valued function of 
$ \mathbf{S} = ( s _{ ij } ) $. 
For a $ p \times p $ Hermitian matrix $ \mathbf{S} = ( s _{ jk} ) $,  
let $ \mathbf{D} _S = ( \partial /  \partial s _{ jk } ) $ be a $ p \times p $ operator matrix, 
the  $ (j,\,k ) $ element of which is given by 
\begin{equation}
\label{eq:2-2a}
  \frac{
    \partial 
  }{
    \partial s _{ jk }
  } 
  =
  \frac{ 1 }{ 2 }
  ( 1 + \delta _{ jk } )
  \left \{
    \frac{ 
      \partial 
    }{
      \partial ( {\rm Re\,} s _{ jk } )
    }
    +
    ( 1 - \delta _{ j k } )
    \sqrt{ -1 }
     \frac{ 
      \partial 
    }{
      \partial ( {\rm Im\,} s _{ jk } )
    }
  \right \}, 
  \qquad
  j,\,k = 1,\,2,\,\ldots,\,p. 
\end{equation}
Here  $ \delta _{ jk } $ is the Kronecker delta ( $=1 \mbox{ if } j = k $ and 
$ =0 \mbox{ if } j \not = k $).  
Thus the $ (j,\,k ) $ element of 
$ \mathbf{D} _S \mathbf{G} ( \mathbf{S} ) $ is  
$$
  \{ \mathbf{D} _S \mathbf{G} ( \mathbf{S} )  \} _{ jk }
  =
  \sum _{ l = 1 } ^p
  \frac{
    \partial g _{ lk  } 
  }{
    \partial s _{ jl }
  }
  ( \mathbf{S} )
  =
   \frac{ 1 }{ 2 }
  ( 1 + \delta _{ jl } )
  \sum _{ l = 1 } ^p
  \left \{
      \frac{ 
      \partial g _{ lk }
    }{
      \partial ( {\rm Re\,} s _{ jl } )
    }
     ( \mathbf{S} )
    +
    ( 1- \delta _{jl} )
    \sqrt{ -1 }
     \frac{ 
      \partial g _{ lk }
    }{
      \partial ( {\rm Im\,} s _{ jl } )
    }
     ( \mathbf{S} )
  \right \}. 
$$
It is directly checked that
$
  \partial s _{ k \ell } /\partial s _{ ij } 
  =
  \delta _{ i \ell } \delta _{ j k },
$ 
and that   
$
  \partial \bar s _{ k \ell } / \partial  s _{ ij } 
  =
  \delta _{ i k } \delta _{ j \ell }
$. 
\vskip 16pt\noindent
\begin{lem}\label{stein-haff}
Assume that 
each entry of $ \mathbf{G} ( \mathbf{S} ) $ is a partially differentiable function with 
respect to $ {\rm Re\,} s _{ j k } $ and $ {\rm Im\,} s _{ j k } $, 
$ j,\,k = 1,\,2,\,\ldots,\,p$.  
Under conditions on $ \mathbf{G}(\mathbf{S}) $ specified in 
{\rm Konno~\cite{Konno2007b}}, the following identity holds: 
\begin{eqnarray}
  \mathbb{E} [
  {\rm Tr } (
    \mathbf{G} (\mathbf{S} )
     \mbox{\boldmath$ \Sigma $} ^{ - 1 }
  )
  ]
  =
  \mathbb{E} [
    ( n - p ) 
     {\rm Tr } (
    \mathbf{G} (\mathbf{S} )
     \mathbf{S} ^{ - 1 }
  )
  +
  {\rm Tr} (
     \mathbf{D} _S \mathbf{G} ( \mathbf{S} )
  )
  ]. 
  \label{eq:2-3}
\end{eqnarray}
\end{lem}

\begin{remark}
{\rm 
The Stein-Haff identity was extended to an elliptically contoured complex distribution by Konno~\cite{Konno2007b}. 
Hence, if we know the improved estimators for the normal case, we can establish the robustness of improvement for  the elliptically contoured complex distribution in a manner similar to that  demonstrated in Kubokawa and Srivastava~\cite{KS1999,KS2001}.  
} 
\end{remark}

\section{Known covariance case}\label{known}
Estimation of a mean matrix  of a real multivariate normal  distribution is considered in 
Stein~\cite{Stein1973}, 
Efron and Morris~\cite{EF1976}, Zheng~\cite{Zheng1986a,Zheng1986b}, and 
Ghosh and Sheih~\cite{GS1991}. Recently Beran~\cite{Beran2007a,Beran2007b} developed adaptive total shrinkage estimators with smaller asymptotic risk than the data matrix. 
\par
In this section, we consider a complex analogue of this problem since the known covariance case gives an insight into estimation problem of  the mean matrix with unknown covariance matrix. 
The problem treated in this section is stated as follows: 
Assume that  $ m \ge p $ and that we observe an $ m \times p $ random matrix 
$ \mathbf{Z} $ with the coordinates $ z _{ ij }\, ( i = 1,\,2,\,\ldots,\,m,\,j = 1,\,2,\,\ldots,\,p )$ that are independently and identically distributed as 
$ \BB{C} N ( \xi _{ ij },\,1 )\,(\xi _{ij } \in \BB{C})$. 
Set $ \mathbf{\Xi} = ( \xi _{ i j } )$, i.e., the $(i,\,j) $ element of 
an $ m \times p $ complex matrix $ \mathbf{\Xi} $ is given by $ \xi _{ij} $. 
We use notation $  \mathbf{Z}\,\sim\,\BB{C} N _{ m \times p } ( 
\mathbf{\Xi},\, \mathbf{I} _m \otimes \mathbf{I} _p ) $ to indicate that 
a random matrix $ \mathbf{Z} $ has a multivariate complex normal distribution with 
a mean matrix $ \mathbf{\Xi} $ and a covariance matrix $ \mathbf{I} _m \otimes \mathbf{I} _p$. 
We consider the problem of estimating the mean matrix $ \mathbf{\Xi} $ under a loss function
\[
  {\SC L } _0( \widehat{\mathbf{\Xi}},\,\mathbf{\Xi} )
  =
  {\rm Tr} \{ ( \widehat{\mathbf{\Xi}} - \mathbf{\Xi} )^*( \widehat{\mathbf{\Xi}} - \mathbf{\Xi} ) \},
\]
where 
$ \widehat{\mathbf{\Xi}} $ is an estimator of $ \mathbf{\Xi} $ based on 
$ \mathbf{Z} $.  
The risk function is given by 
\[
  {\SC R } _0( \widehat{\mathbf{\Xi}},\,\mathbf{\Xi} )
  =
  \BB{E} [ 
    {\rm Tr} \{ ( \widehat{\mathbf{\Xi}} - \mathbf{\Xi} )^*
    ( \widehat{\mathbf{\Xi}} - \mathbf{\Xi} ) \}
  ],
\]
where the expectation is taken with respect to the distribution 
$ \BB{C} N _{ m \times p } ( 
\mathbf{\Xi},\, \mathbf{I} _m \otimes \mathbf{I} _p ) $. 
\subsection{Unbiased risk estimate for a class of invariant estimators}

The maximum likelihood estimator of $ \mathbf{\Xi} $ is given by 
$ \widehat{\mathbf{\Xi}} _{ mle} = \mathbf{Z} $ whose risk function is given by 
$  {\SC R } _0 ( \widehat{\mathbf{\Xi}}_{ mle},\, \mathbf{\Xi} ) = mp $. However, it is expected that the estimator $ \widehat{\mathbf{\Xi}}_{ mle} $ is improved by so-called shrinkage estimators. In order to  search for  
shrinkage estimators in a systematic way, we introduce the following class of estimators and obtain an unbiased risk estimate for this class. This unbiased risk estimate enables us to find a variety of improved estimators.  
\par
Let $ \mathbf{W} = \mathbf{Z} ^* \mathbf{Z} $ and decompose 
$ \mathbf{W} = \mathbf{U} \mathbf{L} \mathbf{U} ^* $, where 
$ \mathbf{U} $ is a $ p \times p $ unitary matrix such that 
$ \mathbf{U} \mathbf{U} ^* = \mathbf{I} _p $ and 
$ \mathbf{L} = {\rm Diag}( \ell _1,\,\ell _2,\,\ldots,\,\ell _p ) $, 
a diagonal real matrix whose $i$-th element $(i = 1,\,2,\,\ldots,\,p)$ is given by 
$ \ell _i$ in the decreasing order. 
Note that all $ \ell _i$'s are non-zero with probability one. 
We consider a class of estimators of the form 
\begin{equation}
  \widehat{\mathbf{\Xi}} _H
  :=
  \widehat{\mathbf{\Xi}} _H ( \mathbf{Z} ) 
  =
  \mathbf{Z} [ \mathbf{I} _p + \mathbf{U} \mathbf{H} ( \mathbf{L} ) \mathbf{U} ^* ], 
\label{eq:2-1}
\end{equation}
where 
$ 
  \mathbf{H}
  := 
  \mathbf{H} (\mathbf{L} ) 
  = 
  {\rm Diag}\, (
     h _1 ( \mathbf{L} ),\,h _2 ( \mathbf{L} ),\,\ldots,\,h _p ( \mathbf{L} ) 
  ) 
  $ with 
$ h _i(\mathbf{L})$'s, $i = 1,\,2,\,\ldots,\,p$, being real-valued functions on $ \BB{R} ^p _{ > }$. 
This class is a complex version of a class of estimators appeared in 
Stein~\cite{Stein1973} and Zheng~\cite{Zheng1986a,Zheng1986b}. 
The following lemma is the complex counterpart of an unbiased risk estimate for orthogonally invariant class of estimators of a mean matrix  of the multivariate real normal distribution, which was proved by Stein~\cite{Stein1973}. 
\begin{lem}\label{steinmatrix}
Assume that 
$ \mathbf{Z}\,\sim\,\BB{C} N _{ m \times p } ( \mathbf{\Xi},\,\mathbf{I} _m \otimes \mathbf{I} _p )$.  
For the estimator $  \widehat{\mathbf{\Xi}} _H $ given by $(\ref{eq:2-1})$, 
we have
\[
  {\SC R } _0( \widehat{\mathbf{\Xi}} _H,\,\mathbf{\Xi} )
  =
  m p 
  + 
  \BB{E} \biggl [
   \sum _{ k = 1 } ^p
   \biggl  \{
     2 ( m - p + 1 ) h _k ( \mathbf{L} ) + 2 \ell _k 
    h _{ kk }
     +
     4 \sum _{ b > k }
     \frac{ 
       \ell _k h _k ( \mathbf{L} )- \ell _b h _b ( \mathbf{L} )
     }{
       \ell _k - \ell _b
     }
     +
     \ell _k h _k ^2( \mathbf{L} )
   \biggr \} 
  \biggr ],
\]
where $ h _{kk} = (\partial h _k / \partial \ell _k )( \mathbf{L} )
\,(k = 1,\,2,\ldots,\,p)$. 
\end{lem}
\par\noindent
{\bf Proof}. 
Using Lemma~\ref{stein} and (\ref{eq:2-1}) we have 
\begin{eqnarray*}
  {\SC R } _0( \widehat{\mathbf{\Xi}} _H,\,\mathbf{\Xi} )
  &=&
  \BB{E} [
    {\rm Tr}\,
    \{
      ( \mathbf{Z} - \mathbf{\Xi} ) ^*( \mathbf{Z} - \mathbf{\Xi} ) 
      +
      2\,{\rm Re}\,(\nabla _Z ^\prime  \mathbf{Z} \mathbf{U} \mathbf{H} \mathbf{U} ^*) +  \mathbf{Z} ^* \mathbf{U} \mathbf{H} ^2 \mathbf{U} ^*  \mathbf{Z} 
    \}
  ]
  \nonumber
  \\
  &=&
  \BB{E} \biggl [
    mp + 
    2\,{\rm Tr}\,
    \{
      {\rm Re}\,(\nabla _Z ^\prime  \mathbf{Z} \mathbf{U} \mathbf{H} \mathbf{U} ^* )\} 
      + \sum _{ i = 1 } ^p \ell _i h _i ^2 
  \biggr ]. 
\end{eqnarray*}
Use Lemma~\ref{lem:a-1a} in the Appendix to evaluate the second term inside expectation of 
the right hand side of the above equation. 
\hfill $\Box$
\begin{remark}
{\rm 
We consider the real version of estimating the mean matrix of the multivariate normal distributions. 
Let  $ \mathbf{X} \sim N _{ m \times p } ( \mathbf{\Xi},\,\mathbf{I} _m \otimes \mathbf{I} _p ) $ 
and decompose $ \mathbf{X} ^\prime \mathbf{X} = \mathbf{O} \mathbf{L} \mathbf{O} ^\prime $ where 
$ \mathbf{O} $ is a $ p \times p$ orthogonal matrix and $ \mathbf{L} = {\rm Diag}( \ell _1,\,\ell,\,\ldots,\,\ell _p ) $ are the ordered eigenvalues of 
$  \mathbf{X} ^\prime \mathbf{X} $ in decreasing order. 
Then the unbiased risk estimate for estimators 
$
  \widehat{\mathbf{\Xi}} _H =  \mathbf{X}[\mathbf{I}_p + \mathbf{O} \mathbf{H} ( \mathbf{L} ) \mathbf{O} ^\prime) ]
$ 
is given by 
\begin{eqnarray*}
  {\SC R } _0( \widehat{\mathbf{\Xi}} _H,\,\mathbf{\Xi} )
  &:=&
  \BB{ E}  [ {\rm Tr}\,( \mathbf{X} - \mathbf{\Xi} ) ^\prime ( \mathbf{X} - \mathbf{\Xi} )]
  \\
  &=&
  \BB{E} \biggl [
      2 ( m - p + 1 ) h _k ( \mathbf{L} ) + 4 \ell _k 
    h _{ kk }
     +
     4 \sum _{ b > k }
     \frac{ 
       \ell _k h _k ( \mathbf{L} )- \ell _b h _b ( \mathbf{L} )
     }{
       \ell _k - \ell _b
     }
     +
     \ell _k h _k ^2( \mathbf{L} )
  \biggr ],
\end{eqnarray*}
which  can be obtained by 
replacing the factor $2$ of $ h _k (\partial h _k / \partial \ell_k )$ in Lemma~\ref{steinmatrix} with 
$4$. 
}
\end{remark}

\subsection{Alternative estimators}
The following proposition is a complex analogue of the results of 
Zheng~\cite{Zheng1986a}. 

\begin{pro}\label{zheng}
Assume that $ m > p $ and 
let $ \gamma _1 ( \mathbf{L} ),\,\gamma _2 ( \mathbf{L} ),\,\ldots,\,\gamma _p ( \mathbf{L} ) $ be functions satisfying 
\par\noindent
{\rm (i)} $ 0 \le \gamma _k ( \mathbf{L} ) \le 2 ( m-p )$; 
\par\noindent
{\rm (ii)} $ (\partial \gamma _k / \partial \ell _k)( \mathbf{L} ) \ge 0 $ for 
$ k = 1,\,2,\,\ldots,\,p$; 
\par\noindent
{\rm (iii)}  $ \gamma _1 ( \mathbf{L} ) \ge \gamma _2 ( \mathbf{L} ) \ge \cdots \ge \gamma _p ( \mathbf{L} ) $. 
\par\noindent
Then the estimator $(\ref{eq:2-1})$ with 
\[
 \mathbf{H}:= 
 \mathbf{H} ( \mathbf{L} ) 
 = 
 -
 {\rm Diag} 
 \left ( 
   \frac{ \gamma _1  ( \mathbf{L} )
   }{
     \ell _1
   }
   ,\,
   \frac{ 
     \gamma _2  ( \mathbf{L} )
   }{
      \ell _2
    }
    ,\,
    \ldots
    ,\,
    \frac{ 
      \gamma _p  ( \mathbf{L} )
    }{ 
      \ell _p
    }
  \right )
 \] 
 is minimax. 
\end{pro}
\par\noindent
{\bf Proof}. 
From Assumptions (i)--(iii) and Lemma~\ref{steinmatrix}, it is easy  
to show that 
$
  \SC{R} _0 ( 
    \mathbf{Z} [ \mathbf{I} _p + \mathbf{U} \mathbf{H} \mathbf{U} ^* ]
    ,\,
    \mathbf{\Xi}
  )
  \le 
  mp  
$.   
\hfill $\Box$
\begin{remark}
{\rm 
Assume that $ m > p $. From Proposition~$\ref{zheng}$, it is easily seen that  
a complex analogue of the crude 
Efron-Morris estimator 
$
   \mathbf{Z} [ \mathbf{I} _p - (m-p) (\mathbf{Z} ^* \mathbf{Z} ) ^{-1} ]
$ 
is minimax. 
}
\end{remark}

\begin{pro}\label{zheng2}
Assume that $ m > p $ and let $ h _k ( \mathbf{L} ) = -(m + p - 2 k ) / \ell _k\,(k = 1,\,2,\,\ldots,\,p )$ in $(\ref{eq:2-1})$. 
Then the estimator of the form $(\ref{eq:2-1})$ is minimax. 
\end{pro}
\par\noindent
{\bf Proof}. 
Without loss of generality we can assume that $ \ell _1 > \ell _2 > \cdots > \ell _p > 0$.  
Let $ h _k ( \mathbf{L} ) = - c _k / \ell _k \,(k = 1,\,2,\,\ldots,\,p )$ in $(\ref{eq:2-1})$, where $ c _k$'s are positive constants such that $ c _1 \ge c _2 \ge \cdots \ge c _p$.  
Then using Lemma~\ref{steinmatrix} and the fact that 
$ \ell _k / ( \ell _k - \ell _b ) > 1 $ for $ b > k $, we can see that 
the risk difference between 
$ \mathbf{Z} $ and $ \mathbf{Z} [ \mathbf{I} _p + \mathbf{U} \mathbf{H} \mathbf{U} ^*] $ is evaluated as 
\begin{eqnarray*}
  \Upsilon 
  &=&
  \SC{R} _0
  (
    \mathbf{Z},\,\mathbf{\Xi}
  )
  -
  \SC{R} _0
  (
    \mathbf{Z} [ \mathbf{I} _p + \mathbf{U} \mathbf{H} \mathbf{U} ^*] ,\,\mathbf{\Xi}
  )
  \\
  &=&
  \sum _{ k = 1 } ^p
  \BB{E} 
  \biggl [
    2 ( m - p ) 
    \frac{ c _k}{ \ell _k } 
    +
     4 \sum _{ b > k }
     \frac{ 
       c _k - c _b 
     }{
       \ell _k - \ell _b
     }
     -
     \frac{ c _k ^2 }{ \ell _k }
  \biggr ]
  \ge
  \sum _{ k = 1 } ^p \BB{E} 
  \left [
    \frac{
      w _k ( c _k )
    }{
      \ell _k
    }
  \right ],
\end{eqnarray*}
where 
\[
  w _k ( t )
  =
  2 ( m + p - 2 k ) t - t ^2 - 4 \sum _{ b = k + 1 } ^p c _b.
\]
If $ c _b = m + p - 2b\,(b = k + 1,\,\ldots,\,p )$, then 
each $ w _k ( t ) $ is maximized at 
$ t = ( m + p - 2 k ) $. Hence, for $ c _k = m + p - 2 b\,(k = 1,\,2,\,\ldots,\,p )$, we can see that 
$ w _k ( m + p - 2 k ) = w _{ k - 1 } ( m + p - 2 k ) < w _{ k - 1 } ( m + p - 2 ( k-1) ) $ for $ k = 2,\,\ldots,\,p $.  Therefore we have 
\[
  0 < w _p ( m - p ) <  w _{ p -1 } ( m - p + 2 ) < \cdots <  w _2 ( m + p - 4 ) < w _1 ( m + p - 2 ),  
\]
from which it follows that $ \Upsilon > 0 $. 
\hfill $\Box$

\begin{remark}\label{transform}
{\rm 
It is easy to extend the result to a known correlated covariance case. Assume that 
we observe an $ m \times p $ random matrix $ \widetilde{\mathbf{Z}} $ 
that is distributed as $ \mathbb{C} N _{ m \times p }( \widetilde{\mathbf{\Xi}},\,\mathbf{I} _m \otimes \mathbf{\Sigma} ) $ with an $ m \times p $ unknown complex matrix $ \widetilde{\mathbf{\Xi}} $ and a known $ p \times p $ positive definite Hermitian matrix $ \mathbf{\Sigma}$. 
Consider the problem of estimating $ \widetilde{\mathbf{\Xi}} $ under the loss function 
$ {\rm Tr} \{ ( \widehat{\widetilde{\mathbf{\Xi}}}- \widetilde{\mathbf{\Xi}} ) ^* ( \widehat{\widetilde{\mathbf{\Xi}}} - \widetilde{\mathbf{\Xi}} ) \mathbf{\Sigma} ^{-1 } \}
$, where 
$ \widehat{\widetilde{\mathbf{\Xi}}} $ is an estimator of 
$ \widetilde{\mathbf{\Xi}}$.  
Transforming 
$ \widetilde{\mathbf{Z}} \to \widetilde{\mathbf{Z}}\mathbf{\Sigma} ^{-1/2}=:\mathbf{Z} $, 
$ \widetilde{\mathbf{\Xi}} \to \widetilde{\mathbf{\Xi}}\mathbf{\Sigma} ^{-1/2}=:\mathbf{\Xi} $, and 
$
\widehat{\widetilde{\mathbf{\Xi}}} \to \widehat{\widetilde{\mathbf{\Xi}}}\mathbf{\Sigma} ^{-1/2}=:\widehat{\mathbf{\Xi}} 
$, the problem reduces to the case when $ \mathbf{\Sigma} = \mathbf{I} _p $. 
Therefore, the Efron-Morris estimator of $ \widetilde{\mathbf{\Xi}} $ is 
given by 
$ \widetilde{\mathbf{Z}} [ \mathbf{I} _p - ( m- p ) (\widetilde{\mathbf{Z}}^*\widetilde{\mathbf{Z}})^{-1} \mathbf{\Sigma} ] $. 
If $\mathbf{\Sigma} $ is unknown and if we observe 
$ \mathbf{S} \sim \BB{C} W _p ( \mathbf{\Sigma},\,n ) $, we replace 
$ \mathbf{\Sigma} $ with $ \mathbf{S} / n$ to obtain an estimator 
$ 
  \widetilde{\mathbf{Z}} [ \mathbf{I} _p - (( m- p )/n) (\widetilde{\mathbf{Z}}^*\widetilde{\mathbf{Z}})^{-1} \mathbf{S} ] 
$. This form of estimators is developed in Section~\ref{unknown}. 
Similarly it is easily seen that the Efron-Morris estimator 
of $ \widetilde{\mathbf{\Xi}} $   is given by 
$ [ \mathbf{I} _m - ( p- m ) (\widetilde{\mathbf{Z}} \mathbf{\Sigma}^{-1} \widetilde{\mathbf{Z}}^*)^{-1}  ] \widetilde{\mathbf{Z}} $ if $ p > m $. 
If $\mathbf{\Sigma} $ is unknown and if we observe 
$ \mathbf{S} \sim \BB{C} W _p ( \mathbf{\Sigma},\,n ) $, we replace 
$ \mathbf{\Sigma} $ with $ \mathbf{S} / n$ to obtain an estimator 
$ 
  [ \mathbf{I} _m - (( p- m )/n) (\widetilde{\mathbf{Z}} \mathbf{S} ^{-1} \widetilde{\mathbf{Z}}^* )^{-1} ] \widetilde{\mathbf{Z}}  
$. This form of estimators is also developed in Section~\ref{unknown}. 
}
\end{remark}

\section{Unknown case and invariant loss}\label{unknown}
In the sequel of this section and the Appendix, we assume 
that $ \mathbf{K} = \mathbf{I} _m $ since the results for the known matrix 
$ \mathbf{K} $ can be obtained from those for $ \mathbf{K} = \mathbf{I} _m $ by 
using a  transformation similar to that given in Remark~\ref{transform}. 
\subsection{Unbiased risk estimate for a class of invariant estimators}

Consider a class of estimators of the form 
$
  \mathbf{Z}
  +
  \mathbf{G} ( \mathbf{Z},\,\mathbf{S} )
$, 
where 
$ \mathbf{G}:=  \mathbf{G} ( \mathbf{Z},\,\mathbf{S} ) $ is an  
$ m \times p $ matrix whose $(i,\,j)$ element $ g _{ ij }\,(i = 1,\,2,\,\ldots,\,m,\,j = 1,\,2,\ldots,\,p )$ is a complex-valued function based on 
 $( \mathbf{Z},\,\mathbf{S} )$.  
\begin{lem}\label{bilodeau-kariya}
Assume that all elements of $ \mathbf{G} ( \mathbf{Z},\,\mathbf{S} ) $ are absolutely continuous functions of $ \mathbf{Z} $ and $ \mathbf{S} $. 
Then we have 
\begin{eqnarray}
  {\SC R} ( 
    \mathbf{Z} + \mathbf{G} ( \mathbf{Z},\,\mathbf{S} )
    ,\,
    (
      \mathbf{\Xi}
      ,\,
      \mathbf{\Sigma}
    )
  )
  &=&
  m p 
  +
  \BB{E} \bigl [
    2\,{\rm Tr\,}\{
      {\rm Re\,} (
        \nabla _Z ^\prime 
        \mathbf{G} ( \mathbf{Z},\,\mathbf{S} )
      )
    \}
    +
    {\rm Tr\,}
    \{
      \mathbf{D} _S 
      \mathbf{G} ^* ( \mathbf{Z},\,\mathbf{S} )
      \mathbf{G} ( \mathbf{Z},\,\mathbf{S} )
    \}
    \nonumber
    \\
    &&
    +
    ( n - p )
    {\rm Tr\,}\{
      \mathbf{G} ^* ( \mathbf{Z},\,\mathbf{S} )
      \mathbf{G} ( \mathbf{Z},\,\mathbf{S} )
      \mathbf{S} ^{-1}
    \} 
  \bigr ]. 
  \label{BK}
\end{eqnarray}
\end{lem}
\par\noindent
{\bf Proof}. 
Use Lemmas~$\ref{stein}$ and $\ref{stein-haff}$.                  
\hfill $\Box$
\par
To describe our class of estimators, 
let $ \mathbf{F} = {\rm Diag\,}( f _1,\,f _2,\,\ldots,\,f _{\min(m,\,p) } ) $ 
be 
the eigenvalues of 
$ \mathbf{Z} ^* \mathbf{Z} \mathbf{S} ^{-1} $.
For $ p > m $  decompose 
$ \mathbf{Z} \mathbf{S} ^{-1} \mathbf{Z}^* = \mathbf{U} \mathbf{F} \mathbf{U} ^* $, where 
$ \mathbf{U} $ is an $ m \times m $ unitary matrix.  For $ m > p $ we decompose 
$ \mathbf{S} = (\mathbf{A}^* )^{-1} \mathbf{A} ^{-1} $ and 
$ 
  \mathbf{Z} ^* \mathbf{Z}  = (\mathbf{A} ^*) ^{ -1 } \mathbf{F} \mathbf{A} ^{-1} $, where $ \mathbf{A} $ is a $ p \times p$ non-singular matrix. 
We consider a class of estimators of the form
\begin{equation}
  \widehat{\mathbf{\Xi}} _H 
  :=
  \widehat{\mathbf{\Xi}} _H ( \mathbf{Z},\,\mathbf{S} )
  =
  \left \{
    \begin{array}{ll}
       \mathbf{Z}
       \{
         \mathbf{I} _p + \mathbf{A}  \mathbf{H}  (\mathbf{F} ) \mathbf{A} ^{-1}
       \}   
     &
       \mbox{if } m > p
      \\ 
       \{
         \mathbf{I} _m + \mathbf{U} \mathbf{H} (\mathbf{F} ) \mathbf{U}^*
       \}\mathbf{Z}
       \qquad
     &
       \mbox{if } p > m
    \end{array}
  \right .
  ,
  \label{eq:3-2}
\end{equation}
where 
$ 
  \mathbf{H}
  := 
  \mathbf{H} (\mathbf{F} ) 
  = 
  {\rm Diag}\, (
     h _1 ( \mathbf{F} ),\,h _2 ( \mathbf{F} ),\,\ldots,\,h _{\min(m,\,p)} ( \mathbf{F} ) 
  ) 
  $ whose $i$-th element 
$ h _i:=h _i ( \mathbf{F} ) $, $i = 1,\,2,\,\ldots,\,\min(m,\,p)$, is a real-valued function on $ \BB{R}  ^{\min(m,\,p)} _{ > } $. 
\par
Let 
\begin{eqnarray}
  \widehat \Delta (n,\,m,\,p;\, \mathbf{H})
  &=&
  \sum _{ k = 1 } ^p
    \biggl \{
      2 ( m - p + 1 ) h _k ( \mathbf{F} )
      +
      2 f _k 
    h _{kk} ( \mathbf{F} ) 
      +
      4 \sum _{ b > k }
      \frac{
        f _k h _k ( \mathbf{F} ) - f _b h _b ( \mathbf{F} )
      }{
        f _k - f _b
      }
      \nonumber
          \\
      &&
      \hskip -12pt
      +
      ( n + p -2 )f _k h _k ^2 ( \mathbf{F} )      -
      2 f _k ^2 h _{ kk } ( \mathbf{F} ) h _k ( \mathbf{F} )
      -
      2 \sum _{ b > k }
      \frac{
        f _k ^2 h _k ^2 ( \mathbf{F} ) - f _b ^2 h _b ^2 ( \mathbf{F} )
      }{
        f _k - f _b
      }
    \biggr \},
    \label{eq:3-2a}
\end{eqnarray}
where  
$ 
  h _{kk} ( \mathbf{F} ) = ( \partial h _k / \partial f _k ) ( \mathbf{F} )
  ,\,  
  k = 1,\,2,\,\ldots,\,p, 
$. 
\begin{pro}\label{pro:main}
Under the suitable conditions, 
we have  
\begin{eqnarray*}
  \SC{R} ( \widehat{\mathbf{\Xi}} _H,\,(\mathbf{\Xi},\,\mathbf{\Sigma} ) )
  &=&
  \left \{
    \begin{array}{ll}
      mp + \BB{E} \bigl [    
        \widehat \Delta (n,\,m,\,p;\, \mathbf{H})
      \bigr ]
      &
      \mbox{if $ m > p$}
    \\
      mp + \BB{E} \bigl [
        \widehat \Delta (n+m-p,\,p,\,m;\, \mathbf{H})
     \bigr ]
      \qquad\qquad
      &
      \mbox{if $ p > m$}
    \end{array}
  \right .
  . 
\end{eqnarray*}
\end{pro}
\par\noindent
{\bf Proof}. 
We apply Lemmas~$\ref{lem:a-3}$ and $\ref{lem:a-5}$ in the Appendix to 
$(\ref{BK})$. 
For $ m > p $, set $ \mathbf{\Phi} = \mathbf{H} $ in the first equation of Lemma~$\ref{lem:a-3}$ and $ \mathbf{\Phi} = \mathbf{F} \mathbf{H} ^2 $ in the second equation of Lemma~$\ref{lem:a-3}$ to get the risk representation for the case when $ m > p$. 
For $ p > m  $, set $ \mathbf{\Phi} = \mathbf{H} $ in the first equation of Lemma~$\ref{lem:a-5}$ and $ \mathbf{\Phi} =  \mathbf{H} ^2 $ in the second equation of Lemma~$\ref{lem:a-5}$ to conclude the proof of the risk representation for the case when $ p > m$. 
\hfill $\Box$

\begin{remark}
{\rm 
Assume that $ \mathbf{\Xi} = \mathbf{0} $ in $(\ref{eq:0})$.
From {\rm \cite{Khatri1965}}, the joint distribution of the eigenvalues of $ \mathbf{Z} ^* \mathbf{Z} \mathbf{S} ^{-1} $ is, aparting from normalizing constants, 
\[
 \Pi _{ k = 1 } ^p
 \frac{
    f _k ^{ m-p} 
 }{
   (1 + f _k) ^{ n + m } 
 }
 \Pi _{ k = 1 } ^{ p-1 } \Pi _{ j = k + 1 } ^p ( f _k - f _j ) ^2
 \Pi _{ k = 1 } ^p {\rm d } f _k
\]
if $ m > p$ while 
it is 
\[
 \Pi _{ k = 1 } ^m
 \frac{
    f _k ^{ p-m} 
 }{
   (1 + f _k) ^{ n+ m  } 
 }
 \Pi _{ k = 1 } ^{ m-1 } \Pi _{ j = k + 1 } ^m ( f _k - f _j ) ^2
 \Pi _{ k = 1 } ^m{\rm d } f _k
\]
if $ p > m $. Note that the substitution rule to get the second distribution from the first distribution, i.e., 
\[
  (p,\,m,\,n) \to (m,\,p,\,n+m-p)
\]
is valid to obtain the second assertion of Proposition~$\ref{pro:main}$ from the first assertion of  Proposition~$\ref{pro:main}$. 
Hence, if we know the estimator of the form 
$  
\mathbf{Z}
       \{
         \mathbf{I} _p + \mathbf{A}  \mathbf{H} \mathbf{A} ^{-1}
       \}   
$ 
when $ m > p $, 
we can easily write down estimators of the form 
$
  \{
         \mathbf{I} _m + \mathbf{U} \mathbf{H} \mathbf{U}^*
       \}\mathbf{Z}
$
when  $ p > m $ by using the above substitution rule.  
}
\end{remark}
\begin{remark}
{\rm 
We consider the real version of estimating the mean matrix of the multivariate normal distributions. 
Let  $ \mathbf{X} \sim N _{ m \times p }  ( \mathbf{\Xi},\,\mathbf{I} _m \otimes \mathbf{\Sigma} _p ) $ 
and $ \mathbf{S} \sim W _p ( n,\,\mathbf{\Sigma})$, where 
$ \mathbf{X} $ and $ \mathbf{S} $ are independent. 
 Let $ \mathbf{F} = {\rm Diag\,}( f _1,\,f _2,\,\ldots,\,f _{\min(m,\,p) } ) $ 
be 
the eigenvalues of 
$ \mathbf{X} ^\prime  \mathbf{X} \mathbf{S} ^{-1} $.
For $ p > m $  decompose 
$ \mathbf{X} \mathbf{S} ^{-1} \mathbf{X}^\prime  = \mathbf{O} \mathbf{F} \mathbf{O} ^\prime $, where 
$ \mathbf{O} $ is an $ m \times m $ orthogonal matrix.  For $ m > p $ we decompose 
$ \mathbf{S} = (\mathbf{A}^\prime  )^{-1} \mathbf{A} ^{-1} $ and 
$ 
  \mathbf{X} ^\prime \mathbf{X}  = (\mathbf{A} ^\prime) ^{ -1 } \mathbf{F} \mathbf{A} ^{-1} $, where $ \mathbf{A} $ is a $ p \times p$ non-singular matrix. 
We consider a class of estimators of the form
$
  \widehat{\mathbf{\Xi}} _H 
  :=
  \mathbf{X}
       \{
         \mathbf{I} _p + \mathbf{A}  \mathbf{H}  (\mathbf{F} ) 
         \mathbf{A} ^{-1}
       \}   
$
where $ m > p $, and 
$ 
  \mathbf{H}
  := 
  \mathbf{H} (\mathbf{F} ) 
  = 
  {\rm Diag}\, (
     h _1 ( \mathbf{F} ),\,h _2 ( \mathbf{F} ),\,\ldots,\,h _{\min(m,\,p)} ( \mathbf{F} ) 
  ) 
  $ whose $i$-th element 
$ h _i:=h _i ( \mathbf{F} ) $, $i = 1,\,2,\,\ldots,\,\min(m,\,p)$, is a real-valued function on $ \BB{R}  ^{\min(m,\,p)} _{ > } $. 
Then the real version of unbiased estimate for the class of estimators 
$\widehat{\mathbf{\Xi}} _H$ is obtained by changing the coefficents of terms 
$ f _k h _{ kk },\,f _k h _k ^2$, and $ f _k h _{kk} h _k $.  
The real version of $ \widehat{\Delta}(n,\,m,\,p;\, \mathbf{H} ) $ in $(\ref{eq:3-2a})$ is 
given as 
\begin{eqnarray*}
  &&
\sum _{ k = 1 } ^p
    \biggl \{
      2 ( m - p + 1 ) h _k ( \mathbf{F} )
      +
      4 f _k 
    h _{kk} ( \mathbf{F} ) 
      +
      4 \sum _{ b > k }
      \frac{
        f _k h _k ( \mathbf{F} ) - f _b h _b ( \mathbf{F} )
      }{
        f _k - f _b
      }
      +
      ( n + p -3 )f _k h _k ^2 ( \mathbf{F} )      
    \nonumber
    \\
    &&
    \qquad
    -
      4 f _k ^2 h _{ kk } ( \mathbf{F} ) h _k ( \mathbf{F} )
      -
      2 \sum _{ b > k }
      \frac{
        f _k ^2 h _k ^2 ( \mathbf{F} ) - f _b ^2 h _b ^2 ( \mathbf{F} )
      }{
        f _k - f _b
      }
      \biggr \}.
\end{eqnarray*} 
}
\end{remark}


\subsection{Alternative estimators}


\begin{pro}\label{zidek}
Let $ \gamma _1 ( \mathbf{F} ),\,\gamma _2 ( \mathbf{F} ),\,\ldots,\,\gamma _{\min(m,\,p)} ( \mathbf{F} ) $ be functions satisfying 
\par\noindent
{\rm (i)} $ 0 \le \gamma _k ( \mathbf{F} ) \le 
\max \{ 2 (m-p)/(n+p),\,2 (p-m)/(n+2m -p ) \}$; 
\par\noindent
{\rm (ii)} $ (\partial \gamma _k / \partial f _k) ( \mathbf{F} )  \ge 0 $ for 
$ k = 1,\,2,\,\ldots,\,\min(m,\,p)$; 
\par\noindent
{\rm (iii)}  $ \gamma _1 ( \mathbf{F} ) \ge \gamma _2 ( \mathbf{F} ) \ge \cdots \ge \gamma _{\min(m,\,p)} ( \mathbf{F} ) $. 
\par\noindent
Then the estimator $(\ref{eq:3-2})$ with 
\[
   \mathbf{H} ( \mathbf{F} ) 
   = 
   -
   {\rm Diag} 
   \left ( 
     \frac{ 
       \gamma _1 ( \mathbf{F} ) 
     }{
       f _1
     }
     ,\,
     \frac{
       \gamma _2 ( \mathbf{F} ) 
     }{
       f _2
     }
     ,\,
     \ldots
     ,\,
     \frac{ 
       \gamma _{ \min(m,\,p ) } ( \mathbf{F} ) 
     }{
        f _{  \min(m,\,p ) } 
      }
    \right )
 \] 
 is minimax. 
\end{pro}
\par\noindent
{\bf Proof}. 
From Assumptions (i)--(iii) and Proposition~\ref{steinmatrix}, it is easy  
to check that 
$
  \SC{R} ( 
    \widehat{\mathbf{\Xi}} _H 
    ,\,
    (\mathbf{\Xi},\,\mathbf{\Sigma})
  )
  \le 
  mp  
$.  
 
\hfill $\Box$

\begin{cor}
The  {\rm Efron-Morris}  estimator 
\[
  \widehat{\mathbf{\Xi}} ^{(EM)}
  =
  \left \{
    \begin{array}{ll}
      \mathbf{Z}
       \{
         \mathbf{I} _p  - \dfrac{m-p}{n+p}   (\mathbf{Z} ^*\mathbf{Z} )^{-1}\mathbf{S}
       \}   
     &
       \mbox{if } m > p
      \\ 
       \{
         \mathbf{I} _m - \dfrac{p-m}{n+2m-p} ( \mathbf{Z} \mathbf{S} ^{-1} \mathbf{Z} ^*)^{-1}
       \}\mathbf{Z}
       \qquad
     &
       \mbox{if } p > m
    \end{array}
  \right .
\] 
is minimax. 
\end{cor}
\par\noindent
{\bf Proof}. 
It is immediately seen from Proposition~\ref{zidek}. 
\hfill $\Box$
\begin{pro}
For $ k = 1,\,2,\,\ldots,\,\min(m,\,p) $, let  
\[
   c ^{(AS) }_k  = \frac{ m + p - 2k }{ n - p + 2 k },
   \qquad
   \mathbf{H} ^{(AS)} ( \mathbf{F} ) 
   =
   -
   {\rm Diag} 
   \left ( 
     \frac{c ^{(AS) }_1}{ f _1 }
     ,\,
     \frac{c ^{(AS) }_2}{ f _2 }
     ,\,
     \ldots
     ,\,
     \frac{c ^{(AS) }_{\min(m,\,p)}
       }{ 
         f _{\min(m,\,p)} 
       }
   \right ).  
\]
Then the estimator
\[
\widehat{\mathbf{\Xi}} ^{(AS)} 
  =
  \left \{
    \begin{array}{ll}
       \mathbf{Z}
       \{
         \mathbf{I} _p + \mathbf{A}  \mathbf{H} ^{(AS)} ( \mathbf{F} )  \mathbf{A} ^{-1}
       \}   
     &
       \mbox{if } m > p
      \\ 
       \{
         \mathbf{I} _m + \mathbf{U} \mathbf{H} ^{(AS)} ( \mathbf{F} )  
         \mathbf{U}^*
       \}\mathbf{Z}
       \qquad
     &
       \mbox{if } p > m
    \end{array}
  \right .
\]
is minimax. 
\end{pro}
\par\noindent
{\bf Proof}. 
It suffices to prove that 
$ 
  \widehat \Delta 
  :=
  \widehat \Delta (n,\,m,\,p;\, \mathbf{H})
  \le 0
$ 
for the case when $ f _1 > f _2 > \cdots > f _p > 0 $ and $ m > p$. 
Put $ h _k = -c _k / f _k $ for $ k = 1,\,2,\,\ldots,\,p $ in $(\ref{eq:3-2a})$, where 
$ c _k$'s are positive constants such that 
$ c _1 \ge c _2 \ge \cdots \ge c _p$. 
Note that 
\begin{eqnarray*}
  \sum _{ k = 1 } \sum _{ b > k }
  \biggl \{
    \frac{ 
      4 (c _k  - c _b)
    }{
      f _k - f _b 
    }
    +
      \frac{ 
      2 (c _k ^2 - c _b ^2)
    }{
      f _k - f _b 
    }
  \biggr \}
  &=&
   \sum _{ k = 1 } \sum _{ b > k }
  \biggl \{
    \frac{ 
      2 (c _k  - c _b)(2 + c _k + c _b )
    }{
      f _k - f _b 
    }
  \biggr \}
  \\
  &\ge&
  \sum _{ k = 1 } \sum _{ b > k }
  \biggl \{
    \frac{ 
       2 (c _k  - c _b)(2 + c _k + c _b )
    }{
      f _k 
    }
  \biggr \}
  \\
  &=&
  \sum _{ k = 1 } 
  \frac{ 1 }{ f _k }
  \biggl \{
    4 ( p - k ) c _k + 2 ( p - k ) c _k ^2
    -
    \sum _{ b > k } \{
      4 c _b + 2 c _b ^2
    \}
  \biggr \}.
\end{eqnarray*}
Hence, for  $ \mathbf{H} = - {\rm Diag} ( c _1/f _1,\,c _2/f _2,\,\ldots,\,c _p /f _p ) $, 
we have 
\begin{eqnarray*}
 \widehat \Delta 
  &\le&
  - 
  \sum _{ k = 1 } ^p
    \frac{ 1 }{ f _k }
    \biggl \{
      2 ( m + p -2k ) c _k
      -
      ( n - p +2k )c _k ^2 
      -
      2 \sum _{ b > k }
      \{
        2 c _b + c _b ^2
      \}
    \biggr \}.
\end{eqnarray*}
Proceed in a  way similar to the proof of Proposition~\ref{zheng2} to see 
that the right hand side of the above inequality is negative. 
\hfill $\Box$


%
%
\vskip 18pt\noindent
{\bf\Large Acknowledgements}
\par\noindent 
This work was in part supported by  the Japan Society for the 
Promotion of Science through Grants-in-Aid for Scientific Research (C) (No.21500283).


\newpage
\appendix
\section{Appendix}
This section develops somewhat tedious computations on eigenvalues, which is a complex analogue of the results given by Loh~\cite{Loh1988,Loh1991a,Loh1991b} and Konno~\cite{Konno1992}. 
In Section~\ref{a1}, we give some results which is needed to prove 
Lemma~\ref{steinmatrix}. In Sections~\ref{a2} and~\ref{a3}, 
we provide some results which is needed to prove Proposition~\ref{pro:main}.  
\par 
In the sequel of this section, we use the following notation:  
For a $ p \times q $ real matrix $ \mathbf{X} = ( x _{ ij } ) _{ i = 1,\,2,
\ldots,\,m,\,j = 1,\,\,2,\,\ldots,\,p} $, define the matrix of differential as 
$ ( {\rm d} \mathbf{X} ) = ( {\rm d } x _{ ij } ) $. We also use the notation 
$ ( {\rm d}\mathbf{X} ) _{ ij }:={\rm d} x _{ij } $.  
For  a $ m \times p $ complex matrix 
$ \mathbf{X} = \mathbf{X} _1 + \sqrt{-1} \mathbf{X} _2 $, 
where $ \mathbf{X} _1,\,\mathbf{X} _2 $ are real matrices, 
we write 
$ ({\rm d} \mathbf{X} ) = ({\rm d} \mathbf{X}_1 )+\sqrt{-1}({\rm d} \mathbf{X}_2 )$.  
\subsection{
Eigencalculus for known covariance case
}~\label{a1}
\hskip -4pt
Let $ \mathbf{W} = \mathbf{Z} ^*\mathbf{Z} = \mathbf{U} \mathbf{L} \mathbf{U} ^* $, where 
$ \mathbf{U} = ( u _{ ij } ) _{ i,\,j = 1,\,\ldots,\,p } $ is a $ p \times p $ unitary matrix and 
$ \mathbf{L} = {\rm Diag} ( \ell _1,\,\ell _2,\,\ldots,\,\ell _p )$ with diagonal elements $ \ell _1,\,\ell _2,\,\ldots,\,\ell _p\,( \ell _1 > \ell _2 > \cdots > \ell _p > 0 )$.  
Recall that \[
  \frac{ 
    \partial 
  }{
    \partial z _{ ij }
  }
  =
  \frac{ 1 }{ 2 }
  \left (
    \frac{ 
      \partial 
    }{
      \partial ({\rm Re\,} z _{ ij } )
    } 
    -
    \sqrt{ -1 }\,
    \frac{ 
      \partial 
    }{
      \partial ({\rm Im\,} z _{ ij } )
    }
  \right ).   
\]
\begin{lem}\label{lem:a-1}
We have 
\begin{eqnarray*}
  \frac{ 
    \partial u _{ i l } 
  }{
    \partial z _{ j k } 
  }
  &=&
  \sum _b \sum _{c \not = l } 
  \frac{ u _{ ic } \bar u _{ bc }  u _{ k l } \bar z _{ j b } }{ \ell _l - \ell _c }, 
  \\
  \frac{ 
    \partial  \bar u _{ i l } 
  }{
    \partial z _{ j k } 
  }
  &=&
  \sum _b \sum _{c \not = l } 
  \frac{ \bar u _{ ic }  \bar u _{ bl }u _{ k c }  \bar z _{ j b }  }{ \ell _l - \ell _c }, 
  \\
  \frac{ 
    \partial \ell _i 
  }{
    \partial z _{ j k } 
  }
  &=&
  \sum _b  u _{ k i } \bar u _{ b i } \bar z _{ j b }. 
\end{eqnarray*}
\end{lem}
\par\noindent
{\bf Proof}. 
Taking the differential of 
$
  \mathbf{W} = \mathbf{U} \mathbf{L} \mathbf{U} ^* 
$ 
we obtain that 
\[
   {\rm d} \mathbf{W} 
   = 
   ({\rm d} \mathbf{U}) \mathbf{L} \mathbf{U} ^* 
   +
    \mathbf{U} \mathbf{L} ({\rm d} \mathbf{U} ^* )
   +
    \mathbf{U} ( {\rm d} \mathbf{L}) \mathbf{U} ^*. 
\]
Multiplying on the left by $ \mathbf{U} ^* $ and on the right by $ \mathbf{U} $ 
we have 
\begin{equation}
  \mathbf{U} ^*
  ( {\rm d} \mathbf{W} )
  \mathbf{U}
  =
  ( \mathbf{U} ^* ( {\rm d} \mathbf{U} ) ) \mathbf{L}
  +
  \mathbf{L} ( \mathbf{U} ^* ( {\rm d} \mathbf{U} ) ) ^*
  +
  {\rm d} \mathbf{L}. 
\label{eq:a-1}
\end{equation}
But, taking the differential of 
$
  \mathbf{U} ^* \mathbf{U} = \mathbf{I} _p 
$, 
we get 
\begin{equation}
  ({\rm d} \mathbf{U} ^*) \mathbf{U} + \mathbf{U} ^* ( {\rm d} \mathbf{U} ) = \mathbf{0}.
\label{eq:a-2}
\end{equation}
Reverting to the coordinates, we obtain from $(\ref{eq:a-1})$ and $(\ref{eq:a-2})$ 
\begin{equation}
  ( \mathbf{U} ^* ( {\rm d} \mathbf{U} ) ) _{ ij }
  =
  \left \{
    \begin{array}{ll}
       \dfrac{ 1 }{
         \ell _j - \ell _i 
       }
       ( \mathbf{U} ^* ( {\rm d} \mathbf{W} ) \mathbf{U} ) _{ij }
       \qquad
       &
       \mbox{for } i \not = j,
       \\
       0 & \mbox{for } i  = j,
    \end{array}
  \right .
\label{eq:a-3}
\end{equation} 
and
\begin{equation}
  ({\rm d} \mathbf{L}) _{ii}
  =
  ( \mathbf{U} ^* ({\rm d} \mathbf{W} ) \mathbf{U} ) _{ ii }. 
  \label{eq:a-4}
\end{equation}
We note that 
\begin{eqnarray}
  ({\rm d} \mathbf{U} ) _{ i l }
  &=&
  \sum _{ c }
  u _{ ic } ( \mathbf{U} ^* ( {\rm d} \mathbf{U}) ) _{ c l }
  =
  \sum _{ c \not = l }
  \frac{
    u _{ ic }
  }{
    \ell _l - \ell _c
  }
  ( \mathbf{U} ^* ( {\rm d} \mathbf{W} ) \mathbf{U} ) _{ c l }
  =
  \sum _{ b _1,\,b _2}
   \sum _{ c \not = l }
  \frac{
    u _{ ic } \bar u _{ b _1 c}  u _{ b _2 l }
  }{
    \ell _l - \ell _c
  }
  ({\rm d} \mathbf{W} ) _{ b _1b _2 }.  
  \nonumber
\end{eqnarray}
But, from 
$
  \mathbf{W} = \mathbf{Z} ^*\mathbf{Z} 
$, 
we observe that 
\begin{eqnarray}
  ({\rm d} \mathbf{W} ) _{ b _1 b  _2 }
  \left (
    \frac{   
      \partial
    }{
      \partial z _{ j k } 
    }
  \right )
  &=&
  \sum _c
  \left \{
    \bar z _{ c b _1 } ( {\rm d} \mathbf{Z} ) _{ c b _2} 
     \left (
    \frac{   
      \partial
    }{
      \partial z _{ j k } 
    }
  \right )
  +
    z _{ c b _2 }
  ({\rm d} \overline{\mathbf{Z}} ) _{ c b _1} 
   \left (
    \frac{   
      \partial
    }{
      \partial z _{ j k } 
    }
  \right )
  \right \}
  =
  \sum _c
  \bar z _{ cb _1 } \delta _{ j c } \delta _{ kb _2}
  \nonumber
  \\
  &=&
  \bar z _{ j b _1 } \delta _{ k b _2},
  \label{eq:a-5}
\end{eqnarray}
from which it follows that 
\begin{eqnarray*}
  \frac{
    \partial u _{ i l }
  }{
    \partial z _{ jk }
  }
  &=&
  ( {\rm d} \mathbf{U} ) _{ i l } 
  \left (
    \frac{ 
      \partial 
    }{
      \partial z _{ jk }
    }
  \right )
  \\
  &=&
  \sum _{ b _1,\,b _2}
  \sum _{ c \not = l }
  \frac{ 
    u _{ i c }   \bar u _{ b _1 c }   u _{ b _2 l }
  }{
    \ell _l - \ell _c
  }
  ( {\rm d} \mathbf{W} ) _{ b _1b _2} 
  \left (
    \frac{
      \partial
    }{
      \partial z _{ jk }
    }
  \right )
  \\
  &=&
  \sum _{ b _1,\,b _2}
   \sum _{ c \not = l }
  \frac{ 
    u _{ i c } \bar u _{ b _1 c } 
     u _{ b _2 l } 
    \bar z _{ j  b _1 } \delta _{ k b _2 }
  }{
    \ell _l - \ell _c
  } 
  \\
  &=&
  \sum _{ b _1 }
  \sum _{ c \not = l }
  \frac{ 
    u _{ i c }  \bar u _{ b _1 c } 
   u _{ k l } \bar z _{ j b _1 }
  }{
    \ell _l - \ell _c
  }. 
\end{eqnarray*}
This completes the first equation of this lemma. 
\par
To prove the second equation we take the complex conjugate of 
$(\ref{eq:a-3})$ and get 
\[
 ( \mathbf{U} ^\prime ( {\rm d} \bar{\mathbf{U}} ) ) _{ ij }
  =
  \left \{
    \begin{array}{ll}
       \dfrac{ 1 }{
         \ell _j - \ell _i 
       }
       ( \mathbf{U} ^\prime ( {\rm d} \overline{\mathbf{W}} ) 
       \overline{\mathbf{U}} ) _{ij }
       \qquad
       &
       \mbox{for } i \not = j,
       \\
       0 & \mbox{for } i  = j.
    \end{array}
  \right .
\]
Using the above equation and noting that  $ \overline{\mathbf{W}} = \mathbf{W} ^\prime $ since 
$ \mathbf{W} $ is Hermitian, we have  
\begin{eqnarray*}
  \frac{
    \partial \bar u _{ i l }
  }{
    \partial z _{ j k }
  }
  &=&
  ({\rm d} \overline{\mathbf{U}} ) _{ i l }
  \left (
    \frac{
      \partial 
    }{
      \partial z _{ j k }
    }
  \right )
  =
  \sum _{ c }
  \bar u _{ ic } ( \mathbf{U} ^\prime ( {\rm d} \overline{\mathbf{U}} ) ) _{ cl }
    \left (
    \frac{
      \partial 
    }{
      \partial z _{ j k }
    }
  \right )
  \\
  &=&
    \sum _{ c \not = l }
     \frac{  \bar u _{ ic }  }{
         \ell _l - \ell _c 
       }
       ( \mathbf{U} ^\prime ( {\rm d} \overline{\mathbf{W}} ) 
       \overline{\mathbf{U}} ) _{cl }
    \left (
    \frac{
      \partial 
    }{
      \partial z _{ j k }
    }
  \right )
  \\
  &=&
  \sum _{ b _1,\,b _2 }
  \sum _{ c \not = l }
     \frac{  
       \bar u _{ ic }   u _{ b _1 c }  \bar u _{ b _2 l }
     }{
         \ell _l - \ell _c 
       }
    ({\rm d} \mathbf{W}) _{ b _2 b _1 }
       \left (
    \frac{
      \partial 
    }{
      \partial z _{ j k }
    }
  \right )   
  \\
  &=&
  \sum _{ b _1,\,b _2}
  \sum _{ c \not = l }
     \frac{  
       \bar u _{ ic }     u _{  b _1 c }   \bar u _{ b _2 l } 
        \bar z _{ j b _2 } \delta _{ k b _1 }
     }{
         \ell _l - \ell _c 
       } 
     \\
  &=&
  \sum _{ b _2}
  \sum _{ c \not = l }
     \frac{ 
       \bar u _{ ic }   
      \bar u _{ b _2 l }u _{ kc }\bar z _{ j b _2} 
     }{
         \ell _l - \ell _c 
       },
\end{eqnarray*}
which completes the proof of the second equation of this lemma. 
The third equality follows from the fact that $ \mathbf{W} $ is Hermitian while the forth equality follows from $(\ref{eq:a-5})$. 
\par
Finally, by $(\ref{eq:a-4})$ and $(\ref{eq:a-5})$, we have  
\begin{eqnarray*}
  \frac{
    \partial \ell _i
  }{
    \partial z _{ j k }
  }
  &=&
  ({\rm d} \mathbf{L} ) _{ii }
  \left (
       \frac{
      \partial 
    }{
      \partial z _{ j k }
    }
  \right )
  =
   ( \mathbf{U} ^* ({\rm d} \mathbf{W} ) \mathbf{U} ) _{ ii } 
   \left (
       \frac{
      \partial 
    }{
      \partial z _{ j k }
    }
  \right )
  \\
  &=&
  \sum _{ b _1,\,b _2 }
  \bar u _{ b _1 i }   u _{ b _2 i } 
  ({\rm d} \mathbf{W} ) _{ b _1 b _2 }
    \left (
       \frac{
      \partial 
    }{
      \partial z _{ j k }
    }
  \right )
  \\
  &=&
  \sum _{ b _1,\,b _2 }
  \bar u _{ b _1 i } u _{ b _2 i} \bar z _{ j b _1 } \delta _{k b _2 } 
  =
  \sum _{ b _1 }
  \bar u _{ b _1 i }  u _{ ki } \bar z _{ j b _1 },
\end{eqnarray*}
which completes the proof of the third equation of this lemma.   
\hfill $\Box$

\begin{lem}\label{lem:a-1a}
Let 
$ 
  \mathbf{\Phi} ( \mathbf{L} )
  =
  {\rm Diag\,}
  ( 
    \varphi _1 ( \mathbf{L} )
    ,\,
    \varphi _2 ( \mathbf{L} )
    ,\,
    \ldots
    ,\,
     \varphi _p ( \mathbf{L} )
   )
$, 
where 
$    \varphi _i ( \mathbf{L} )$'s $( i = 1,\,2,\,\ldots,\,p )$ are 
differentiable functions from $ \BB{R} ^p _{ > } \to \BB{R} _+ $. 
Then we have 
\begin{eqnarray*}
  {\rm Tr\,} \{ {\rm Re}\,(
    \nabla _Z ^\prime \mathbf{Z}  \mathbf{U} \mathbf{\Phi} (\mathbf{L}) \mathbf{U} ^*
  )\}
  &=&
  \sum _k 
  \left \{
    (m-p+ 1) \varphi _k (\mathbf{L})
    +
    2 \sum _{ c > k }
    \frac{ \ell _k \varphi _k (\mathbf{L}) - \ell _c \varphi _c (\mathbf{L}) }{ \ell _k - \ell _c}
    +
    \ell _k 
    \frac{ 
       \partial \varphi _k
    }{
      \partial \ell _k
    }
    (\mathbf{L})
  \right \}. 
\end{eqnarray*}
\end{lem}

\par\noindent
{\bf Proof}. 
Write $ \mathbf{\Phi} $ and $ \varphi _i $ for $ \mathbf{\Phi}(\mathbf{L} ) $ and  $ \varphi _i (\mathbf{L} )\,(i = 1,\,2,\,\ldots,\,p) $, respectively.  
Note that  
\begin{equation}
  {\rm Tr\,}\{{\rm Re\,}  ( \nabla _Z ^\prime \mathbf{Z} \mathbf{U} \mathbf{\Phi} \mathbf{U} ^* )\}
  =
  m {\rm Tr\,} \mathbf{H} + 
  \frac{ 1 }{ 2 }
  \biggl ({\rm Tr\,} ( \mathbf{Z} ^\prime \nabla _Z (\mathbf{U} \mathbf{\Phi} \mathbf{U} ^*)^\prime )
  +
  {\rm Tr\,} ( \mathbf{Z} ^* \nabla _{\bar Z} \mathbf{U} \mathbf{\Phi} \mathbf{U} ^* )
  \biggr ),  
  \label{eq:2-2b}
\end{equation}
where $ \nabla _{\bar Z } = ( \partial / \partial { \bar z } _{ j k } ) _{ j = 1,\,2,\,\ldots,\,m,\,k = 1,\,2,\ldots,\,p } $ with 
$ \partial / \partial \bar z _{ j k } = \{ \partial / \partial ({\rm Re\,} z _{ jk } )
+ \sqrt{-1} \partial / \partial ({\rm Im\,} z _{ jk })\} / 2$.  
We use Lemma~$\ref{lem:a-1}$ to evaluate 
the second term in the right hand side of $(\ref{eq:2-2b})$ as  
\begin{eqnarray*}
  {\rm Tr\,} ( 
    \mathbf{Z} ^\prime \nabla _{ Z} (\mathbf{U} \mathbf{\Phi} \mathbf{U} ^* 
  ) ^\prime )
  &=&
  \sum _{ i,\,j,\,k,\,l}
  z _{ ji }
  \frac{ \partial (\bar u _{ kl } \varphi _l u _{ i l })
  }{
    \partial z _{ j k }
  }
  \\
  &=&
  \sum _{i,\,j,\,k,\,l}
  z _{ ji }
  \biggl \{
    \varphi _l u _{ i l }
    \frac{ 
      \partial  \bar u _{ k l }
    }{
      \partial z _{ j k }
    }
    +
    \varphi _l \bar u _{ k l }
    \frac{ 
      \partial u _{i l }
    }{
      \partial z _{ j k }
    }
    +
    \bar u _{ k l } u _{ i l }
    \sum _{ k ^\prime}  
    \frac{ \partial \varphi _l }{ \partial \ell _{ k ^\prime} }
    \frac{ 
      \partial \ell _{ k ^\prime} 
    }{
      \partial z _{ j k }
    }
  \biggr \}
  \\
  &=&
  \sum _{i,\,j,\,k,\,l}
  z _{ ji }
  \biggl \{
  \varphi _l u _{ i l } 
  \sum _b \sum _{c \not = l }
  \frac{
    \bar u _{ k c }  u _{ k c } \bar u _{ b l } \bar z _{ j b } 
  }{
    \ell _l - \ell _c
  }
  +
  \varphi _l \bar u _{ k l }
  \sum _b \sum _{c \not = l }
  \frac{ 
    u _{ i c }  \bar u _{ b c }
  u _{ kl}  \bar z _{ jb }
  }{
    \ell _l - \ell _c
  } 
  \\
  &&
  \hskip 84pt
  +
  \bar u _{ k l } u _{ i l }
  \sum _{{ k ^\prime} ,\,b }
  u _{ k{ k ^\prime}  } 
  \bar u _{ b { k ^\prime}  } \bar z _{ jb }
  \frac{
    \partial \varphi _l
  }{
    \partial \ell _{ k ^\prime} 
  }
  \biggr \}
  \\
  &=&
  \sum _k 
  \left \{ 
    \sum _{ c \not = k }
    \frac{ \ell _k  \varphi _k - \ell _k  \varphi _c }{ \ell _k - \ell _c}
    +
    \ell _k 
    \frac{ 
       \partial  \varphi _k
    }{
      \partial \ell _k
    }
  \right \}. 
\end{eqnarray*}
We use Lemma~$\ref{lem:a-1}$ and the fact that 
$ 
  \partial q / \partial \bar{z} 
  =
  \overline{
      \partial \bar{q} / \partial  z 
  }
$ to evaluate 
the third term in the right hand side of $(\ref{eq:2-2b})$ as  
\begin{eqnarray*}
  {\rm Tr\,} ( 
    \mathbf{Z} ^* \nabla _{ \bar Z} \mathbf{U} \mathbf{\Phi} \mathbf{U} ^*) )
  &=&
  \sum _{ i,\,j,\,k,\,l}
  \bar z _{ ji }
  \frac{ \partial ( u _{ kl } \varphi _l \bar u _{ i l })
  }{
    \partial \bar z _{ j k }
  }
  \\
  &=&
  \sum _{i,\,j,\,k,\,l}
  \bar z _{ ji }
  \biggl \{
    \varphi _l \bar u _{ i l }
    \overline { 
      \frac{ 
        \partial  \bar u _{ k l }
      }{
        \partial z _{ j k }
      }
    }
    +
    \varphi _l  u _{ k l }
    \overline{
      \frac{ 
        \partial u _{i l }
      }{
        \partial  z _{ j k }
      }
    }
    +
    u _{ k l } \bar u _{ i l }
    \sum _{ k ^\prime}  
    \frac{ \partial \varphi _l }{ \partial \ell _{ k ^\prime} }
    \overline{
      \frac{ 
        \partial \ell _{ k ^\prime} 
      }{
        \partial  z _{ j k }
      }
    }
  \biggr \}
  \\
  &=&
  \sum _k 
  \left \{ 
    \sum _{ c \not = k }
    \frac{ \ell _k  \varphi _k - \ell _k  \varphi _c }{ \ell _k - \ell _c}
    +
    \ell _k 
    \frac{ 
       \partial  \varphi _k
    }{
      \partial \ell _k
    }
  \right \}. 
\end{eqnarray*}
Putting the above two equations into $(\ref{eq:2-2b})$, we have 
\begin{eqnarray*}
 {\rm Tr\,} \{  {\rm Re\,}(
    \nabla _Z ^\prime \mathbf{Z} \mathbf{U} \mathbf{\Phi} 
    \mathbf{U} ^* 
  )\}
  &=&
  \sum _k 
  \left \{
    m  \varphi _k 
    +
    \sum _{ c \not = k }
    \frac{ \ell _k  \varphi _k - \ell _c  \varphi _c +  \varphi _c ( \ell _c - \ell _k ) }{ \ell _k - \ell _c}
    +
    \ell _k 
    \frac{ 
       \partial  \varphi _k
    }{
      \partial \ell _k
    }
  \right \}  
  \\
  &=&
  \sum _k 
  \left \{
    (m-p+ 1)  \varphi _k 
    +
    \sum _{ c \not = k }
    \frac{ \ell _k  \varphi _k - \ell _c  \varphi _c  }{ \ell _k - \ell _c}
    +
    \ell _k 
    \frac{ 
       \partial  \varphi _k
    }{
      \partial \ell _k
    }
  \right \}.  
\end{eqnarray*}
Combining this equation with $(\ref{eq:2-2b})$, we completes the proof of this lemma. 
\hfill $\Box$

\subsection{
Eigencalculus for unknown covariance case with $ m > p$
}~\label{a2}
Next we record calculus on the eigenvalues for the case when $ m > p $. 
Let $ \mathbf{A} = ( a _{ ij } ) _{ i,\,j = 1,\,2,\ldots,\,p } $ be a 
$ p \times p $ nonsingular matrix such that 
\[
  \mathbf{A} ^* \mathbf{S} \mathbf{A} 
  =
  \mathbf{I} _p
  ,\qquad
  \mathbf{A} ^* \mathbf{Z} ^* \mathbf{Z} \mathbf{A} 
  =
  \mathbf{F}
  ,\qquad
  \mathbf{F}
  =
  {\rm Diag\,} ( f _1,\,f _2,\,\ldots,\,f _p )
\]
with $ f _1 > f _2 > \cdots > f _p > 0$. This means that we consider 
$ (\mathbf{Z},\,\mathbf{S} ) $ such that  a matrix  $ \mathbf{Z} ^* \mathbf{Z} \mathbf{S} ^{-1} $ has the distinct eigenvalues 
$ f _1,\,f _2,\,\ldots,\,f _p $. 

\begin{lem}\label{lem:a-2}
Let $ \mathbf{A} ^{-1} = ( a ^{ ij } ) _{ i,\,j = 1,\,2,\ldots,\,p } $, 
$ \overline{\mathbf{A}}  = ( \bar a _{ ij } ) _{ i,\,j = 1,\,2,\ldots,\,p } $, 
and $ (\overline{\mathbf{A}}) ^{-1} = ( \bar a ^{ ij } ) _{ i,\,j = 1,\,2,\ldots,\,p } $. 
For $ i,\,k,\,k ^\prime = 1,\,2,\,\ldots,\,p$, and $ j = 1,\,2,\,\ldots,\,m $, 
we have 
\begin{eqnarray*}
  \frac{
    \partial a ^{ l k ^\prime }
  }{
    \partial z _{ j k }
  } 
  &=&
  \sum _b \sum _{ c \not = l }
  \frac{
     \bar a _{b l } a _{ k c } a ^{ c k ^\prime } \bar z _{ j b } 
  }{
    f _l - f _c 
  },
  \\
  \frac{
    \partial  a _{ i l }
  }{
    \partial z _{ j k }
  }
  &=&
  \sum _b \sum _{ c \not = l }
  \frac{
     a _{i c } \bar a _{ b c } a _{ kl } \bar z _{ j b } 
  }{
    f _l - f _c 
  },
  \\
  \frac{
    \partial  f _{k ^\prime} 
  }{
    \partial z _{ j k }
  } 
  &=&
  \sum _{ b  }
  \bar a _{b k ^\prime } a _{ k k ^\prime }  \bar z _{ j b },
  \\
   \frac{
    \partial ( a ^{ k i } \bar a ^{kj}  )
  }{
    \partial s _{ ij }
  }  
  &=&
   a ^{ ki } \bar  a ^{ kj}  
  \bar a _{ j k }  a _{ i k } 
  +
  \sum _{ b  \not = k}
  a ^{ ki }  
  a _{ i k }
  \bar a _{ j b } 
  \bar a ^{ b j } 
  \frac{ 
    f _{ b  }
  }{
    f _{ b  } - f _k
  }
   + 
  \sum _{ b  \not = k}
  \bar a ^{ kj} 
  \bar a _{ j k }
  a _{ i b  } 
  a ^{ b  i } 
  \frac{ 
    f _{ b  }
  }{
    f _{ b  } - f _k
  },
  \\
   \frac{
    \partial  f _i
  }{
    \partial s _{ k k ^\prime }
  } 
  &=&
  - \bar a _{ k ^\prime  i } a _{ k i } f _i.
\end{eqnarray*}
\end{lem}
\par\noindent
{\bf Proof}. 
Put 
$
  \mathbf{W} = \mathbf{Z} ^* \mathbf{Z}
  =
  ( w _{ ij } ) _{ i,\,j= 1,\,2,\,\ldots,\,p }
$. 
Differentiating   
$ \mathbf{S}  = (\mathbf{A} ^*) ^{-1} \mathbf{A}^{-1}$ 
and 
$  \mathbf{W} = (\mathbf{A} ^*)^{-1}\mathbf{F}  \mathbf{A}^{-1}$, 
we have 
\begin{eqnarray*}
  ( {\rm d} \mathbf{S} )
  &=&
  (\mathbf{A} ^*) ^{-1} ({\rm d} \mathbf{A}^{-1})
  +
  ( {\rm d} (\mathbf{A} ^*)^{-1} ) \mathbf{A}^{-1}
  \\
  ({\rm d} \mathbf{W})
  &=&
   (\mathbf{A} ^*) ^{-1} \mathbf{F} ({\rm d} \mathbf{A}^{-1})
  +
  ( {\rm d} (\mathbf{A} ^*)^{-1} ) \mathbf{F} \mathbf{A}^{-1}
  +
  (\mathbf{A} ^*)^{-1} ( {\rm d} \mathbf{F})  \mathbf{A}^{-1}.
\end{eqnarray*}
Multiplying these equations by $ \mathbf{A} ^* $ on the left and 
by $ \mathbf{A} $ on the right, we get 
\begin{eqnarray}
   \mathbf{A} ^* ( {\rm d} \mathbf{S} )  \mathbf{A} 
  &=&
 ({\rm d} \mathbf{A}^{-1}) \mathbf{A} 
  +
   \mathbf{A} ^* ( {\rm d} (\mathbf{A} ^*)^{-1} ) 
  \label{eq:a-6}
  \\
  \mathbf{A} ^* ({\rm d} \mathbf{W}) \mathbf{A} 
  &=&
  \mathbf{F} ({\rm d} \mathbf{A}^{-1}) \mathbf{A} 
  +
 \mathbf{A} ^*  ( {\rm d} (\mathbf{A} ^*)^{-1} ) \mathbf{F} 
  +
  ( {\rm d} \mathbf{F}).
   \label{eq:a-7}
\end{eqnarray}
To obtain the derivatives with respect to $ \partial / \partial z _{ jk } $, 
we may assume that $ {\rm d} \mathbf{S} = 0 $.  
Then, putting $(\ref{eq:a-6})$ into $(\ref{eq:a-7})$ and from some algebraic calculation, we have 
\begin{equation}
  ({\rm d} \mathbf{F} ) _{ k ^\prime k ^\prime }
  =
  \sum _{ b,\,c }
  \bar a _{ b k ^\prime } 
  a _{ c k ^\prime}
  ( {\rm d} \mathbf{W} ) _{ b c } 
\label{eq:a-8}
\end{equation}
and
\begin{equation}
  ( ( {\rm d} \mathbf{A} ^{-1} ) \mathbf{A} ) _{ l b }
  =
  \left \{
    \begin{array}{ll}
      \dfrac{
        1
      }{
        f _l - f _b 
      }
      \sum  _{ c _1,\,c _2 }
      \bar a _{ c _1 l }  a _{ c _2 b } ({\rm d} \mathbf{W} ) _{ c _1 c _2 } 
      \qquad
      & 
      \mbox{if }
      l \not = b,
      \\
      0 & \mbox{if }
      l = b.  
    \end{array}
  \right .
  \label{eq:a-9}
\end{equation}
From $(\ref{eq:a-9})$ and the fact that 
\begin{eqnarray*}
  ({\rm d} \mathbf{W} ) _{ c _1 c _2 }
  \left (
    \frac{
      \partial 
    }{
      \partial z _{ j k }
    }
  \right )
  &=&
  \sum _{ c _3 }
  \biggl \{
     z _{ c _3 c _2 }
    ({\rm d} \overline{\mathbf{Z} } )  _{ c _3 c _1 }
    \left (
      \frac{
        \partial 
      }{
        \partial z _{ j k }
      }
    \right )
    +
    \bar z _{ c _3 c _1 }
    ({\rm d} \mathbf{Z} ) _{ c _3 c _2 }
     \left (
      \frac{
        \partial 
      }{
        \partial z _{ j k }
      }
    \right )
  \biggr \}
  \\
  &=&
  \sum _{ c _3 } \bar z _{ c _3 c _1 } \delta _{ c _3 j } \delta _{ c _2 k }
  =
  \bar z _{ j c _1 } \delta _{ c _2 k }, 
\end{eqnarray*}
we have 
\begin{eqnarray*}
  \frac{ 
    \partial  a ^{ l k ^\prime }
  }{
    \partial z _{ j k }
  }  
  &=&
  ( {\rm d} \mathbf{A} ^{-1 } ) _{ l k ^\prime }
    \left (
      \frac{
        \partial 
      }{
        \partial z _{ j k }
      }
    \right )
  =
  \sum _b 
    a ^{ b k ^\prime }
   ( ({\rm d} \mathbf{A} ^{-1 }) \mathbf{A}  ) _{ l b } 
    \left (
      \frac{
        \partial 
      }{
        \partial z _{ j k }
      }
    \right )
  \\
  &=&
  \sum _{ c _1,\,c _2 }
  \sum _{ b \not = l }
  \frac{ 
     \bar a _{ c _1 l }   a _{ c _2 b }
      a ^{ b k ^\prime}
   }{ 
    f _l  - f _b 
  }
  ( {\rm d} \mathbf{W} ) _{ c _1 c _2 } 
     \left (
      \frac{
        \partial 
      }{
        \partial z _{ j k }
      }
  \right )
  \\
  &=&
   \sum _{ c _1 }
    \sum _{ b \not = l } 
  \frac{ 
    \bar a _{ c _1 l }  
    a _{ k b }
    a ^{ b k ^\prime}
    \bar z _{ j c _1 }
  }{ 
    f _l  - f _b 
  }, 
\end{eqnarray*}
which completes the first equation of this lemma. 
\par
Differentiating $ \mathbf{A} \mathbf{A} ^{-1} = \mathbf{I} _p $, we have
$
  ( {\rm d} \mathbf{A}) \mathbf{A} ^{-1}
  +
  \mathbf{A} ( {\rm d} \mathbf{A} ^{-1})
  =
  \mathbf{0}
$. 
Multiplying this equation by $ \mathbf{A} $ on the right
 and 
using the first assertion of this lemma, we have 
\begin{eqnarray*}
  \frac{
    \partial   a _{ i l }
  }{
    \partial z _{ j k }
  }
  &=&
  -
  \sum _{ b _1,\,b _2 }
  a _{ i b _1 }  a _{ b _2 l }
  \frac{ 
    \partial a ^{ b _1 b _2 }
  }{
    \partial z _{ j k }
  } 
  =
  -
  \sum _{ b _1,\,b _2 }
  a _{ i b _1 } a _{ b _2 l } 
  \sum _{ b _4 } \sum _{ b _3 \not = b _1}
  \frac{ 
     \bar a _{ b _4 b _1 } a _{ k b _3 }
     a ^{ b _3 b _2 } \bar z _{ j b _4 }
   }{ f _{ b _1 } - f _{ b _3 } }
  \\
  &=&
  -
  \sum _{ b _4 } \sum _{ b _1 \not = l }
  \frac{
    a _{ i b _1 } \bar a _{ b _4 b _1 }
    a _{ k l } \bar z _{ j b _4 }
  }{ 
    f _{ b _1 } - f _{ l} 
  },
\end{eqnarray*}
which completes the proof of the second assertion. 
\par
From $(\ref{eq:a-8})$ we have 
\begin{eqnarray*}
  \frac{ 
    \partial f _{ k ^\prime }
  }{
    \partial z _{ jk }
  }
  &=&
  \sum _{ b _1,\,b _2 }
  \bar a _{ b _1 k ^\prime }
   a _{ b _2 k ^\prime }
  ({\rm d} \mathbf{W} ) _{ b _1 b _2 }
  \left (
    \frac{
      \partial 
    }{
      \partial z _{ jk }
    }
  \right ) 
  =
  \sum _{ b _1,\,b _2 }
  \bar a _{ b _1 k ^\prime} 
  a _{ b _2 k ^\prime }
  \bar z _{ j b _1 } \delta _{ b _2 k }
  =
  \sum _{ b _1 }
  \bar a _{ b _1 k ^\prime }
  a _{ k k ^\prime }
  \bar z _{ j b _1 },
\end{eqnarray*}
which completes the third assertion of this lemma. 
\par
To derive the derivatives with respect to $ s _{ ij } $ we assume that 
$ {\rm d} \mathbf{W} = 0$ in $(\ref{eq:a-6})$. Reverting to the coordinates, we have 
\begin{equation}
  ({\rm d} \mathbf{F})  _{ii}
  =
  - 
  \sum _{ j,\,k } \bar a _{ ji } a _{ k i } f _i ( {\rm d} \mathbf{S} ) _{ jk } 
  \label{eq:a-10}
\end{equation}
and 
\begin{eqnarray}
  ( ({\rm d} \mathbf{A} ^{-1} ) \mathbf{A} ) _{ ij }
  &=&
  \frac{
    f _j
  }{
    f _j - f _i
  }
  ( \mathbf{A} ^* ( {\rm d} \mathbf{S} ) \mathbf{A} ) _{ ij },
  \qquad
  \mbox{if } i \not = j,
  \label{eq:a-11}
  \\
  ( \overline{({\rm d} \mathbf{A} ^{-1} ) \mathbf{A}} ) _{ ij }
  &=&
  \frac{
    f _j
  }{
    f _j - f _i
  }
  ( \overline{\mathbf{A} ^* ( {\rm d} \mathbf{S} ) \mathbf{A}} ) _{ ij },
  \qquad
  \mbox{if } i \not = j.
  \label{eq:a-12}
\end{eqnarray}
Since 
$
 \mathbf{A} ^* ( {\rm d} \mathbf{S} ) \mathbf{A}
  =
  ( {\rm d} \mathbf{A} ^{-1} ) \mathbf{A} + \mathbf{A} ^* 
  ( {\rm d} (\mathbf{A}^*)^{-1} ) 
$ which implies that 
$
  2 {\rm Re\,} [(( {\rm d} \mathbf{A} ^{-1} ) \mathbf{A} ) _{ii}] 
  =
  ( \mathbf{A} ^* ( {\rm d} \mathbf{S} ) \mathbf{A})_{ii}. 
$,  
we have 
\begin{equation}
  ( ( {\rm d} \mathbf{A} ^{-1} ) \mathbf{A} ) _{ ii }
  +
  ( \overline{( {\rm d} \mathbf{A} ^{-1} ) \mathbf{A}} ) _{ ii }
  =
  ( \mathbf{A} ^* ( {\rm d} \mathbf{S} ) \mathbf{A})_{ii}.  
\label{eq:a-13}
\end{equation}
From $(\ref{eq:a-11})-(\ref{eq:a-13})$ we have 
\begin{eqnarray*}
  \frac{
    \partial ( a ^{ k i } \bar a ^{kj}  )
  }{
    \partial s _{ ij }
  }  
  &=&
  a ^{ki} 
  \frac{
    \partial \bar a ^{kj} 
  }{
    \partial s_{ ij }
  }
  +
  \bar a ^{kj} 
  \frac{
    \partial a ^{ki}
  }{
    \partial s _{ ij }
  } 
    =
   a ^{ki} 
  (\overline{{\rm d} \mathbf{A}^{-1}} ) _{ k j } 
  \left (
    \frac{
      \partial 
    }{
      \partial s_{ ij }
    }
  \right )
  +
  \bar  a ^{kj}  
  ({\rm d} \mathbf{A}^{-1} ) _{ki} 
  \left (
    \frac{
      \partial 
    }{
      \partial s_{ ij }
    }
  \right )
   \\
  &=&
  a ^{ ki } 
  \sum _{ b _1 }
   \bar a ^{ b _1 j }
  (\overline{ ({\rm d} \mathbf{A} ^{-1}) \mathbf{A} }) _{ k b _1 } 
  \left (
    \frac{
      \partial 
    }{
      \partial s_{ ij }
    }
  \right ) 
  +
  \bar a ^{ kj} 
    \sum _{ b _1 }
     a ^{ b _1 i } 
   (({\rm d} \mathbf{A} ^{-1}) \mathbf{A} ) _{ k b _1 } 
  \left (
    \frac{
      \partial 
    }{
      \partial s_{ ij }
    }
  \right )
 \\
 &=&
 a ^{ ki } 
 \bar a ^{ kj}  
 \left \{  (\overline{ ({\rm d} \mathbf{A} ^{-1}) \mathbf{A} }) _{ k k } 
     \left (
    \frac{
      \partial 
    }{
      \partial s_{ ij }
    }
  \right )
   +
   (({\rm d} \mathbf{A} ^{-1}) \mathbf{A} ) _{ k k } 
     \left (
    \frac{
      \partial 
    }{
      \partial s_{ ij }
    }
  \right )
  \right \}
 \\
 &&
 +
   a ^{ ki } 
  \sum _{ b _1 \not = k}
   \bar a ^{ b _1 j }
  (\overline{ ({\rm d} \mathbf{A} ^{-1}) \mathbf{A} }) _{ k b _1 } 
  \left (
    \frac{
      \partial 
    }{
      \partial s_{ ij }
    }
  \right ) 
   +
  \bar a ^{ kj} 
    \sum _{ b _1 \not = k}
     a ^{ b _1 i } 
   (({\rm d} \mathbf{A} ^{-1}) \mathbf{A} ) _{ k b _1 } 
  \left (
    \frac{
      \partial 
    }{
      \partial s_{ ij }
    }
  \right )
  \\
  &=&
  a ^{ ki } \bar a ^{ kj}  
  \sum _{ b _2,\,b _3}
  \bar a _{ b _2 k }  a _{ b _3 k } 
  ( {\rm d} \mathbf{S} ) _{ b _2 b _3 }
   \left (
    \frac{
      \partial 
    }{
      \partial s_{ ij }
    }
  \right )
  \\
  &&
  +
  a ^{ ki } 
  \sum _{ b _2,\,b _3}
  \sum _{ b _1 \not = k}
  a _{ b _2 k }
  \bar a _{ b _3 b _1 } 
  \bar  a ^{ b _1 j } 
  \frac{ 
    f _{ b _1 }
  }{
    f _{ b _1 } - f _k
  }
  ({\rm d} \overline{\mathbf{S}} ) _{ b _2 b _3 }
  \left (
    \frac{
      \partial 
    }{
      \partial s_{ ij }
    }
  \right ) 
  \\
  &&
   +
  \bar a ^{ kj} 
  \sum _{ b _2,\,b _3}
   \sum _{ b _1 \not = k}
  \bar a _{ b _2 k }
  a _{ b _3 b _1 } 
  a ^{ b _1 i } 
  \frac{ 
    f _{ b _1 }
  }{
    f _{ b _1 } - f _k
  }
  ({\rm d} \mathbf{S} ) _{ b _2 b _3 }
  \left (
    \frac{
      \partial 
    }{
      \partial s_{ ij }
    }
  \right ) 
  \\
  &=&
  a ^{ ki } \bar a ^{ kj}  
  \sum _{ b _2,\,b _3}
  \bar  a _{ b _2 k }  
  a _{ b _3 k } 
  \delta _{ b _2 j } \delta _{ b _3 i }
  +
  a ^{ ki } 
  \sum _{ b _2,\,b _3}
  \sum _{ b _1 \not = k}
  a _{ b _2 k }
  \bar  a _{ b _3 b _1 } 
  \bar a ^{ b _1 j } 
  \frac{ 
    f _{ b _1 }
  }{
    f _{ b _1 } - f _k
  }
  \delta _{ b _2 i } \delta _{ b _3 j }
  \\
  &&
   +
  \bar  a ^{ kj} 
  \sum _{ b _2,\,b _3}
  \sum _{ b _1 \not = k}
  \bar a _{ b _2 k }
  a _{ b _3 b _1 } 
  a ^{ b _1 i } 
  \frac{ 
    f _{ b _1 }
  }{
    f _{ b _1 } - f _k
  }
  \delta _{ b _2 j } \delta _{ b _3 i } 
    \\
  &=&
  a ^{ ki } \bar  a ^{ kj}  
  \bar a _{ j k }  a _{ i k } 
  +
  \sum _{ b _1 \not = k}
  a ^{ ki }  
  a _{ i k }
  \bar a _{ j b _1 } 
  \bar a ^{ b _1 j } 
  \frac{ 
    f _{ b _1 }
  }{
    f _{ b _1 } - f _k
  }
  \\
  &&
   + 
  \sum _{ b _1 \not = k}
  \bar a ^{ kj} 
  \bar a _{ j k }
  a _{ i b _1 } 
  a ^{ b _1 i } 
  \frac{ 
    f _{ b _1 }
  }{
    f _{ b _1 } - f _k
  },
\end{eqnarray*}
which completes the proof of forth assertion of this lemma.
\par
Finally, from $(\ref{eq:a-10})$, we have
\begin{eqnarray*}
  \frac{
    \partial f _{b _1 }
  }{
    \partial s _{ ij }
  }
  &=&
  ( {\rm d} \mathbf{F} ) _{ b _1 b _1 } 
  \left (
    \frac{
      \partial 
    }{
      \partial s _{ ij }
    }
  \right )
  =
  -
  \sum _{ b _2,\,b _3 }
  \bar a _{ b _2 b _1 } 
  a _{ b _3 b _1 }
  f _{ b _1 }
  ({\rm d} \mathbf{S}) _{ b _2 b _3 }
    \left (
    \frac{
      \partial 
    }{
      \partial s _{ ij }
    }
  \right )
  =
   -
  \sum _{ b _2,\,b _3 }
  \bar a _{ b _2 b _1 } 
  a _{ b _3 b _1 }
  f _{ b _1 }
  \delta _{ b _2 j } 
  \delta _{ b _3 i }
  \\
  &=&
  -
  \bar a _{ j b _1 }  a _{ i b _1 } f _{ b _1 },
\end{eqnarray*}
which completes the final part of this lemma. 
\hfill $\Box$

\begin{lem}\label{lem:a-3}
Let 
$ 
  \mathbf{\Phi} ( \mathbf{F} )
  =
  {\rm Diag\,}
  ( 
    \varphi _1 ( \mathbf{F} )
    ,\,
    \varphi _2 ( \mathbf{F} )
    ,\,
    \ldots
    ,\,
     \varphi _p ( \mathbf{F} )
   )
$, 
where 
$    \varphi _i ( \mathbf{F} )$'s $( i = 1,\,2,\,\ldots,\,p )$ are 
differentiable functions from $ \BB{R} ^p _{ > } \to \BB{R} _+ $. 
Then we have 
\begin{eqnarray*}
  {\rm Tr\,} \{{\rm Re\,}(
    \nabla _Z ^\prime 
     \mathbf{Z} 
      \mathbf{A} 
      \mathbf{\Phi} 
      \mathbf{A} ^{-1})
  \}
  &=&
  \sum _k
  \biggl \{
    f _k \varphi _{kk} ( \mathbf{F} )
    +
     (m -  p + 1 ) \varphi _k( \mathbf{F} )
    +
    2 \sum _{ b > k }
    \frac{
      f _k \varphi _k ( \mathbf{F} )
      -
       f _b \varphi _b ( \mathbf{F} )
    }{
      f _k - f _b
    }
  \biggr \},
  \\
  {\rm Tr\,}
  (
    \mathbf{D} _S (
      (\mathbf{A} ^*) ^{ -1 }
      \mathbf{\Phi} ( \mathbf{F} )
      \mathbf{A} ^{-1} 
    )
  )
  &=&
  \sum _k 
  \biggl \{
    ( 2 p - 1 ) \varphi _k ( \mathbf{F} )
    -
    2 \sum _{ b > k }
    \frac{
      f _k \varphi _k ( \mathbf{F} )
      -
      f _b \varphi _b ( \mathbf{F} )
    }{
      f _k - f _b
    }
    - 
    f _k \varphi _{kk} ( \mathbf{F} )
  \biggr \},
\end{eqnarray*}
where 
$
  \varphi _{ kk } ( \mathbf{F} ) 
  =
  (\partial \varphi _k / \partial f _k)( \mathbf{F} )
  ,\,
  k = 1,\,2,\,\ldots,\,p. 
$
\end{lem}
\par\noindent
{\bf Proof}. 
Use notation  $ \mathbf{\Phi} $, $ \varphi _k $, and  
$ \varphi _{kk} $ short for 
$ \mathbf{\Phi} ( \mathbf{F} )$, $ \varphi _{ k } ( \mathbf{F} )  $, and $
  \varphi _{ kk } ( \mathbf{F} ) $, respectively.  
To prove the first equation of this lemma, 
we first note that
\[
  {\rm Tr\,} \{{\rm Re\,}(
    \nabla _Z ^\prime 
     \mathbf{Z} 
      \mathbf{A} 
      \mathbf{\Phi} 
      \mathbf{A} ^{-1})
  \}
  =
  m {\rm Tr\,} \mathbf{\Phi} 
  +
  \frac{ 1 }{ 2 }
  \biggl (
  {\rm Tr\,} \{
    \mathbf{Z} ^\prime
    \nabla _Z (
      (\mathbf{A} ^\prime)^{-1}
      \mathbf{\Phi} 
      \mathbf{A} ^\prime 
    )
    \}
    +
    {\rm Tr\,} \{
    \mathbf{Z} ^*
    \nabla _{\bar Z} (
      (\mathbf{A} ^*)^{-1}
      \mathbf{\Phi} 
      \mathbf{A} ^* 
      )
    \}
    \biggr ). 
\]
Now 
we use the first three equations in 
Lemma~$\ref{lem:a-2}$ to evaluate the second term inside expectation of the right hand side as 
\begin{eqnarray*}
  {\rm Tr\,} \{
    \mathbf{Z} ^\prime
    \nabla _Z (
      (\mathbf{A} ^\prime)^{-1}
      \mathbf{\Phi} 
      \mathbf{A} ^\prime 
    )
  \}
  &=&
  \sum _{ i,\,j,\,k,\,l }
  z _{ ji }
  \frac{
    \partial ( a ^{ lk } \varphi _l a _{ i l } )
  }{
    \partial z _{ jk }
  }
  \\
  &=&
  \sum _{ i,\,j,\,k,\,l }
   z _{ ji }
  \left \{ 
    \varphi _l a _{ i l }
    \frac{
      \partial  a ^{ lk } 
    }{
      \partial z _{ jk }
    }
    +
    \varphi  _l a ^{ l k }
    \frac{
      \partial a _{ i l } 
    }{
      \partial z _{ jk }
    }
    +
    a ^{ lk }  a _{ i l }
    \frac{
      \partial  \varphi _l 
    }{
      \partial z _{ jk }
    }
  \right \}
  \\
  &=&
  \sum _{ i,\,j,\,k,\,l }
   \biggl \{
  z _{ ji }
  \varphi _l a _{ i l }
  \sum _{ b _2 } \sum _{ b _1 \not = l}
  \frac{
    \bar a _{ b _2 l }
    a _{ k b _1 }
    a ^{ b _1 k }
    \bar z _{ jb _2 }
  }{
    f _l - f _{ b _1 }
  }
  +
  z _{ ji } \varphi  _l a ^{ l k }
  \sum _{ b _2 } \sum _{ b _1 \not = l}
  \frac{
    a _{ i b _1  }
    \bar a _{ b _2 b _1 }
    a _{ k l }
    \bar z _{ jb _2 }
  }{
    f _l - f _{ b _1 }
  }
    \\
  &&
  \qquad\qquad
  +
  z _{ ji }
  a ^{ lk }  a _{ i l }
  \sum _{ b _1,\,b _2 }
  \bar a _{ b _2 b _1 }
  a _{ k b _1 }
  \bar z _{ j b _2 }
  \frac{
    \partial  \varphi _l 
  }{
    \partial f _{ b _1 }
  }  
  \biggr \}
  \\
  &=&
  \sum _{ i,\,j,\,l,\,b _2 }
  \biggl \{
  \sum _{ b _1 \not = l}
  \bar z _{ jb _2 }
   z _{ ji }a _{ i l }\bar a _{ b _2 l }
  \frac{
  \varphi _l 
  }{
    f _l - f _{ b _1 }
  }
  +
  \sum _{ b _1 \not = l}
   \bar z _{ jb _2 } z _{ ji }
    a _{ i b _1  }
  \bar a _{ b _2 b _1 }
  \frac{
   \varphi  _l 
  }{
    f _l - f _{ b _1 }
  }
    \\
  &&
  \qquad\qquad
  +
  \bar z _{ j b _2 }
  z _{ ji }
   a _{ i l }
  \bar a _{ b _2 l }
  \frac{
    \partial  \varphi _l 
  }{
    \partial f _{ l }
  }  
  \biggr \}
  \\
  &=&
  {\rm Tr\,} (
    \mathbf{Z} ^* \mathbf{Z} 
    \mathbf{A} \widetilde{\mathbf{\Phi}} 
    \mathbf{A} ^*
  ),
\end{eqnarray*}
where 
$  
  \widetilde{\mathbf{\Phi}}
  =
  {\rm Diag\,}
  (
    \widetilde{\mathbf{\phi}} _1     
    ,\,
     \widetilde{\mathbf{\phi}} _2 
     ,\,
     \ldots
     ,\,
      \widetilde{\mathbf{\phi}} _p 
  )
$ 
with 
$
   \widetilde{\mathbf{\phi}} _i 
   =
   (\partial \varphi _i/\partial f _i )
   +
   \sum _{ b \not = i }
   ( \varphi _i  - \varphi _b ) / ( f _i - f _b )
   ,\
   (i = 1,\,2,\,\ldots,\,p )
$. 
Similarly  
we use the first three equations in 
Lemma~$\ref{lem:a-2}$ to evaluate the third term inside expectation of the right hand side as 
\begin{eqnarray*}
  {\rm Tr\,} \{
    \mathbf{Z} ^*
    \nabla _{ \bar Z} (
      (\mathbf{A} ^*)^{-1}
      \mathbf{\Phi} 
      \mathbf{A} ^* 
    )
  \}
  &=&
  \sum _{ i,\,j,\,k,\,l }
  \bar z _{ ji }
  \frac{
    \partial ( \bar a ^{ lk } \varphi _l \bar a _{ i l } )
  }{
    \partial \bar z _{ jk }
  }
  \\
  &=&
  \sum _{ i,\,j,\,k,\,l }
  \bar z _{ ji }
  \left \{ 
    \varphi _l \bar a _{ i l }
    \overline{
      \frac{
        \partial a ^{ lk } 
      }{
      \partial z _{ jk }
      }
    }
    +
    \varphi  _l \bar a ^{ l k }
    \overline{ 
      \frac{
        \partial a _{ i l } 
      }{
        \partial z _{ jk }
      }
    }
    +
    \bar a ^{ lk }  \bar a _{ i l }
    \overline{
      \frac{
        \partial  \varphi _l 
      }{
        \partial z _{ jk }
      }
    }
  \right \}
  \\
  &=&
  {\rm Tr\,} (
    \mathbf{Z} ^* \mathbf{Z} 
    \mathbf{A} \widetilde{\mathbf{\Phi}} 
    \mathbf{A} ^*
  ).
\end{eqnarray*}
Putting $ \mathbf{Z} ^* \mathbf{Z} = (\mathbf{A} ^* )^{-1} \mathbf{F} \mathbf{A} $ into the right hand side of the above two equations, we have 
\begin{eqnarray*}
   {\rm Tr\,} \{{\rm Re\,}(
    \nabla _Z ^\prime 
     \mathbf{Z} 
      \mathbf{A} 
      \mathbf{\Phi} 
      \mathbf{A} ^{-1})
  \}
  &=&
  \sum _k 
  \biggl \{
    f _k \varphi _{ kk }
    +(m-p+1 )
    \varphi _k 
    +
    \sum _{ b \not = k }
    \frac{
      f _k \varphi _k - f _b \varphi _b  
    }{
      f _k - f _b
    }
  \biggr \}, 
\end{eqnarray*}
which completes the first equation of this lemma. 
\par
Next we prove the second equation of this lemma. 
Apply the chain rule first and use the forth and fifth equations 
of Lemma~$\ref{lem:a-2}$ 
to get 
\begin{eqnarray*}
   {\rm Tr\,}
  (
    \mathbf{D} _S (
      (\mathbf{A} ^*) ^{ -1 }
      \mathbf{\Phi} 
      \mathbf{A} ^{-1} 
    )
  )
  &=&
  \sum _{ i,\,j,\,k }
  \frac{
    \partial ( \bar a ^{ kj } \varphi _k a ^{ ki })
  }{
    \partial s _{ ij }
  }
  \\
  &=&
  \sum _k
  \biggl \{
    \varphi _k 
    \sum _{ i,\,j }
     \frac{
    \partial ( \bar a ^{ kj }  a ^{ ki })
  }{
    \partial s _{ ij }
  }
  +
  \sum _{ i,\,j }
  a ^{ k i } \bar a ^{ kj } 
  \sum _{ b  }
  \frac{
    \partial \varphi _k
  }{
    \partial f _{ b }
  }
  \frac{
    \partial f _{ b  }
  }{
    \partial s _{ ij } 
  }
  \biggr \}
  \\
  &=&
  \sum _k
  \biggl \{
    \varphi _k
    \biggl (
      1 
      +
     2 \sum _{ b  \not = k }
      \frac{
        f _{ b  }
      }{
        f _{ b } - f _k
      }
    \biggr )
    -
    f _k
     \frac{
    \partial \varphi _k
  }{
    \partial f _{ k }
  }
  \biggr \}
  \\
  &=&
  \sum _k
  \biggl \{
   \varphi _k
   +
   2 \sum _{ b  \not = k }
   \frac{
     ( f _{ b  } - f _k ) \varphi _k + f _k \varphi _k 
   }{
     f _{ b } - f _k
   }
   -
   f _k \varphi _{ kk }
  \biggr \}
    \\
  &=&
  \sum _k
  \biggl \{
   \varphi _k
   +
   2(p-1)\varphi _k
   -
   2 \sum _{ b  > k }
   \frac{
     f _k \varphi _k - f _{ b  } \varphi _{ b  } 
   }{
     f _{ k } - f _{ b  }
   }
   -
   f _k \varphi _{ kk },
  \biggr \},
\end{eqnarray*} 
which completes the proof of this lemma. 
\hfill $\Box$

\subsection{
Eigencalculus for unknown covariance case with $ m < p$ 
}~\label{a3}
Let $ \mathbf{Z} ^* \mathbf{S} ^{-1 } \mathbf{Z} = \mathbf{U} \mathbf{F} \mathbf{U} ^* $, where 
$ \mathbf{U} = ( u _{ ij } ) _{ i,\,j = 1,\,\ldots,\,m } $ is an $ m \times m $ unitary matrix and 
$ \mathbf{F} = {\rm Diag} ( f _1,\,f _2,\,\ldots,\,f _m )$ with diagonal elements $ f _1,\,f _2,\,\ldots,\,f _m\,( f _1 > f _2 > \cdots > f _m > 0 )$.  
\begin{lem}\label{lem:a-4}
For $ i,\,k,\,l,\,l ^{\prime} = 1,\,2,\,\ldots,\,p,\,j = 1,\,2,\,\ldots,\,m $, we have 
\begin{eqnarray*}
    \frac{ 
    \partial u _{ il }
  }{
    \partial z _{ j k  }
  }
  &=&
  \sum _{ b _3,\,b _5 }
   \sum _{ b _1 \not = l }
  \frac{
    u _{ i b _1 } \bar u _{ j b _1 } u _{ b _3 l } 
    s ^{ k b _5 } \bar z _{ b _3 b _5 }
  }{
    f _l - f _{ b _1 }
  },
  \\
    \frac{
    \partial \bar u _{ il }
  }{
    \partial z _{ j k }
  }
  &=&
  \sum _{ b _3,\,b _5 }
   \sum _{ b _1 \not = l }
  \frac{
    \bar u _{ ib _1} \bar u _{ j l }  u _{ b _3 b _1 }
    s ^{ k b _5 } \bar z _{ b _3 b _5}
  }{
      f _l - f _{ b _1 }
  },
  \\
    \frac{ 
    \partial f _{ b _1 }
  }{
    \partial z _{ j k  }
  }
  &=&
  \sum _{ b_3,\,b _5}
  \bar u _{ j b _1 } u _{ b _3 b _1 } 
  s ^{ k b _5} \bar z _{ b _3 b _5},
  \\
  \frac{
    \partial u _{ k  l }
  }{
    \partial s _{ ij }
  }
  &=&
  -
  \sum _{ b _2,\,b _3,\,b _4,\,b _5 }
  \sum _{ b _1 \not = l }
  \frac{
    u _{ k b _1 }
    \bar u _{ b _2 b _1 }
    u _{ b _3 l}
    z _{ b _2 b _4 }
    \bar z _{ b _3 b _5 }
    s ^{ b _4 j }
    s  ^{ i b _5} 
  }{
    f _l - b _{ b _1 }
  },
  \\
  \frac{
    \partial \bar u _{ k  l }
  }{
    \partial s _{ ij }
  }
  &=&
  -
  \sum _{ b _2,\,b _3,\,b _4,\,b _5 }
  \sum _{ b _1 \not = l }
  \frac{
    \bar u _{ k b _1 }
    \bar u _{ b _2 l }
    u _{ b _3  b _1 }
    z _{ b _2 b _4 }
     \bar z _{ b _3 b _5 }
    s ^{ b _4 j }
    s  ^{ i b _5} 
  }{
    f _l - b _{ b _1 }
  },
   \\
  \frac{
    \partial f _{ l ^{\prime} }
  }{
    \partial s _{ ij }
  }
  &=&
  -
  \sum _{ b _2,\,b _3,\,b _4,\,b _5 }
    \bar u _{ b _2 l ^{\prime} }
    u _{ b _3 l ^{\prime}  }
    z _{ b _2 b _4 }
    s ^{ b _4 j }
    s  ^{ i b _5}
    \bar z _{ b _3 b _5 }.
\end{eqnarray*}

\end{lem}
\par\noindent
{\bf Proof}. 
Let 
$
  \mathbf{T}
  =
  \mathbf{Z} \mathbf{S} ^{-1} \mathbf{Z} ^*
$ and take differential of $ \mathbf{T} = \mathbf{U} \mathbf{F} \mathbf{U} ^*$ to get 
\[
  ({\rm d} \mathbf{T})
  =
  ({\rm d} \mathbf{U})
  \mathbf{F} 
  \mathbf{U} ^*
  +
  \mathbf{U} 
  \mathbf{F} 
  ({\rm d} \mathbf{U} ^*)
  +
  \mathbf{U} 
  ({\rm d} \mathbf{F} )
  \mathbf{U} ^*,   
\]
from which it follows that 
\[
  \mathbf{U} ^* ( {\rm d} \mathbf{T} ) \mathbf{U}
  =
  \mathbf{U} ^* ( {\rm d} \mathbf{U} ) \mathbf{F}
  +
  \mathbf{F} ( {\rm d} \mathbf{U} ^* ) \mathbf{U}
  +
  ( {\rm d} \mathbf{F} )
  =
  ( \mathbf{U} ^* ({\rm d} \mathbf{U} ) ) \mathbf{F}
  -
  \mathbf{F} ( \mathbf{U} ^* ({\rm d} \mathbf{U} ) )
  +
  ( {\rm d} \mathbf{F} ). 
\]
Therefore we have 
\begin{equation}
  ( \mathbf{U} ^* ( {\rm d} \mathbf{U} ) ) _{ ij }
  =
  \left \{
    \begin{array}{ll}
      \dfrac{ 1 }{ f _j - f _i }
      ( \mathbf{U} ^* ( {\rm d} \mathbf{T} ) \mathbf{U} ) _{ ij }
      \qquad
      & \mbox{ if } i \not = j
      \\
      0 & \mbox{ if } i = j,
    \end{array}
  \right. 
  \label{eq:a-13a}
\end{equation}
while, for $ i = j $, we have
\begin{eqnarray}
  ({\rm d} \mathbf{F} )  _{ i i }
  =
  ( \mathbf{U} ^*( {\rm d} \mathbf{T} ) \mathbf{U} ) _{ ii }. 
  \label{eq:a-13b}
\end{eqnarray}
Hence 
\begin{eqnarray*}
  \frac{ 
    \partial u _{ il }
  }{
    \partial z _{ j k  }
  }
  &=&
  ({\rm d} \mathbf{U}) _{ il }
  \left (
    \frac{ 
      \partial 
    }{
      \partial z _{ j k  }
    }
  \right )
  =
  \sum _{ b _1}
  u _{ i b _1 }
  ( \mathbf{U} ^* ( {\rm d} \mathbf{U} ) ) _{ b _1 l }
   \left (
    \frac{ 
      \partial 
    }{
      \partial z _{ j k  }
    }
  \right )
  \\
  &=&
  \sum _{ b _2,\,b _2 }
  \sum _{ b _1 \not = l }
  \frac{
    u _{ i b _1 }
    \bar u _{ b _2 b _1 }
  u _{ b _3 l }
  }{
    f _l - f _{ b _1 }
  }
  ( {\rm d} \mathbf{T} ) _{ b _2 b _3 }
  \left (
   \frac{ 
      \partial 
   }{
      \partial z _{ j k  }
   }
   \right ).
\end{eqnarray*}
But 
\begin{eqnarray}
  ( {\rm d} \mathbf{T} ) _{ b _2 b _3 }
  \left (
    \frac{ 
      \partial 
    }{
      \partial z _{ j k  }
    }
   \right )
   &=&
   \sum _{ b _4,\,b _5 }
   \left \{
     \bar z _{ b _3 b _5 }
     s ^{ b _4 b _5 }
     ( {\rm d} \mathbf{Z} ) _{ b _2 b _4 }
     \left (
       \frac{ 
          \partial 
       }{
         \partial z _{ j k  }
       }
     \right )
     +
     z _{ b _2 b _4 }
     s ^{ b _4 b _5 }
     ( {\rm d} \bar{\mathbf{Z}} ) _{ b _3 b _5 }
     \left (
       \frac{ 
         \partial 
       }{
         \partial z _{ j k  }
       }
     \right )
  \right \}
  \nonumber
  \\
  &=&
  \sum _{ b _4,\,b _5 }
  \delta _{ j b _2 } \delta _{ k b _4} 
  s ^{ b _4 b _5 }
  \bar z _{ b _3 b _5 }
  =
  \sum _{ b _5}
  \delta _{ j b _2 } s ^{ k b _5 } \bar z _{ b _3 b _5}.
  \label{eq:a-14}
\end{eqnarray}
Therefore we have 
\begin{eqnarray*}
    \frac{ 
    \partial u _{ il }
  }{
    \partial z _{ j k  }
  }
  &=&
  \sum _{ b _2,\,b _3,\,b _5 }
  \sum _{ b _1 \not = l }
  \frac{
    u _{ i b _1 } \bar u _{ b _2 b _1 } u _{ b _3 l } \delta _{ j b _2 }
    s ^{ k b _5 } \bar z _{ b _3 b _5 }
  }{
    f _l - f _{ b _1 }
  }
  =
  \sum _{ b _3,\,b _5 }
  \sum _{ b _1 \not = l }
  \frac{
    u _{ i b _1 } \bar u _{ j b _1 } u _{ b _3 l } 
    s ^{ k b _5 } \bar z _{ b _3 b _5 }
  }{
    f _l - f _{ b _1 }
  },
\end{eqnarray*}
which complete the proof of the first assertion of this lemma. 
\par
Similarly we have  
$
  ( ( {\rm d} \mathbf{U} ^*) \mathbf{U} ) _{ ij }
  =
  ( \mathbf{U} ^* ( {\rm d} \mathbf{T} ) \mathbf{U} ) _{ ij } / ( f _i - f _j )
$ for $ i \not = j$,  
from which it follows that 
\begin{eqnarray*}
  \frac{
    \partial \bar u _{ il }
  }{
    \partial z _{ j k }
  }
  &=&
  ( {\rm d} \mathbf{U} ^* ) _{ li } 
  \left (
    \frac{
      \partial 
    }{
      \partial z _{ j k }
    }
  \right )
  =
  \sum _{ b _1 } 
  \bar u _{ i b _1 }
  ( ( {\rm d} \mathbf{U} ^*) \mathbf{U} ) _{ l b _1 } 
  \left (
    \frac{
      \partial 
    }{
      \partial z _{ j k }
    }
  \right )
  = 
  \sum _{ b _2,\,b _3 }
   \sum _{ b _1 \not = l }
  \frac{\bar u _{ i b _1 }
  \bar u _{ b _2l}
  u _{ b _3 b _1 }}{ f _l - f _{ b _1 }}
  ( {\rm d} \mathbf{T} ) _{ b _2 b _3 }
  \left (
    \frac{
      \partial 
    }{
      \partial z _{ j k }
    }
  \right )
  \\
  &=&
  \sum _{ b _2,\,b _3,\,b _5 }
  \sum _{ b _1 \not = l }
  \frac{
    \bar u _{ i b _1 } \bar u _{ b _2l}   u _{ b _3 b _1 }
    \delta _{ j b _2 } s ^{ k b _5 } \bar z _{ b _3 b _5}
  }{
     f _l - f _{ b _1 }
  }
     \\
  &=&
  \sum _{ b _3,\,b _5 }
  \sum _{ b _1 \not = l }
  \frac{
    \bar u _{ i b _1 }  \bar u _{ j l }  u _{ b _3 b _1 }
    s ^{ k b _5 } \bar z _{ b _3 b _5}
  }{
      f _l - f _{ b _1 }
  },
\end{eqnarray*}
which completes the second assertion of this lemma. 
\par
Furthermore, from $(\ref{eq:a-13b})$ and $(\ref{eq:a-14})$,  we have  
\begin{eqnarray*}
  \frac{ 
    \partial f _{ b _1 }
  }{
    \partial z _{ j k  }
  }
  &=&
  ( {\rm d} \mathbf{F} ) _{ b _1 b _1 }
  \left (
    \frac{
      \partial 
    }{
      \partial z _{ j k }
    }
  \right )
  =
  ( \mathbf{U} ^* ( {\rm d} \mathbf{T} ) \mathbf{U} ) _{ b _1 b _1 }
  \left (
    \frac{
      \partial 
    }{
      \partial z _{ j k }
    }
  \right )
  =
  \sum _{ b _2,\,b _3 }
  \bar u _{ b _2 b _1 }
  u _{ b _3 b _1 }
  ( {\rm d} \mathbf{T} ) _{ b _2 b _3 }
    \left (
    \frac{
      \partial 
    }{
      \partial z _{ j k }
    }
  \right )
  \\
  &=&
  \sum _{ b _2,\,b_3,\,,b _5}
  \bar u _{ b _2 b _1 } u _{ b _3 b _1 } \delta _{ j b _2 }
  s ^{ k b _5} \bar z _{ b _3 b _5}
  =
  \sum _{ b_3,\,,b _5}
  \bar u _{ j b _1 } u _{ b _3 b _1 } 
  s ^{ k b _5} \bar z _{ b _3 b _5},
\end{eqnarray*}
which complete the proof of the third assertion of this lemma. 
\par
Next we prove the forth equation of this lemma. 
Use $(\ref{eq:a-13a})$ to get 
\begin{eqnarray*}
  \frac{
    \partial u _{ k l }
  }{
    \partial s _{ ij }
  }
  &=&
  \sum _{ b _1 \not = l }
  u _{ k b _1 } 
  ( \mathbf{U} ^* ( {\rm d} \mathbf{U}  ) ) _{ b _1 l }
  \left (
    \frac{
      \partial
    }{
      \partial s _{ ij }
    }
  \right )
  \\
  &=&
   \sum _{ b _2,\,b _3 }
   \sum _{ b _1 \not = l }
  \frac{
     u _{ k b _1 }
     \bar u _{ b _2 b _1 }
    u _{ b _3 l }
  }{
    f _l - f _{ b _1 }
  }
  ({\rm d} \mathbf{T} ) _{ b _2 b _3 }
    \left (
    \frac{
      \partial
    }{
      \partial s _{ ij }
    }
  \right ).
\end{eqnarray*}
But, since 
$
  ( \partial s ^{ b _4 b _5 }  / \partial s _{ ij}) 
  = - s ^{b _4j } s ^{i b _5 }  
$, we have 
\[
   ({\rm d} \mathbf{T} ) _{ b _2 b _3 }
    \left (
    \frac{
      \partial
    }{
      \partial s _{ ij }
    }
  \right )
  =
  \sum _{ b _4,\, b _5 }
  z _{ b _2 b _4 }
  \bar z _{ b _3 b _5 }
  ( {\rm d} \mathbf{S} ^{-1}) _{ b _4 b _5 }
     \left (
    \frac{
      \partial
    }{
      \partial s _{ ij }
    }
  \right )
  =
  -
  \sum _{ b _4,\, b _5 }
  z _{ b _2 b _4 } s ^{ b _4 j } s ^{ i b _5 } 
  \bar z _{ b _3 b _5 }.
\] 
Combining these two equations we completes the proof of the forth equation of this lemma. 
\par
Similarly we have 
\begin{eqnarray*}
  \frac{
    \partial \bar u _{ k l }
  }{
    \partial s _{ ij }
  }
  &=&
  \sum _{ b _1 }
  \bar u _{ k b _1 } 
  ( ( {\rm d} \mathbf{U} ^*) \mathbf{U} ) _{ l b _1 }
  \left (
    \frac{
      \partial
    }{
      \partial s _{ ij }
    }
  \right )
  \\
  &=&
   \sum _{ b _2,\,b _3 }
   \sum _{ b _1 \not = l }
  \frac{
     \bar u _{ k b _1 }
    \bar u _{ b _2 l }
    u _{ b _3 b _1 }
  }{
    f _l - f _{ b _1 }
  }
  ({\rm d} \mathbf{T} ) _{ b _2 b _3 }
    \left (
    \frac{
      \partial
    }{
      \partial s _{ ij }
    }
  \right )
  \\
  &=&
  -
  \sum _{ b _2,\,b _3,\,b _4,\,b _5 }
  \sum _{ b _1 \not = l }
  \frac{
    \bar  u _{ k b _1 }
    \bar u _{ b _2 l }
    u _{ b _3 b _1 }
    z _{ b _2 b _4} 
     \bar z _{ b _3 b _5}
    s ^{ b _4 j } s ^{ i b _5 }
  }{
    f _l - f _{ b _1 }
  },  
\end{eqnarray*}
which completes the proof of the forth assertion of this lemma. 
\par
Finally, from $(\ref{eq:a-13b})$, we have
\begin{eqnarray*}
  \frac{
    \partial f _{ b _1 }
  }{
    \partial s _{ ij }
  }
  &=&
  \sum _{ b _2,\,b _3 }
  \bar u _{ b _2 b _1 } u _{ b _3 b _1 }
  ( {\rm d} \mathbf{T} ) _{ b _2 b _3 }
      \left (
    \frac{
      \partial
    }{
      \partial s _{ ij }
    }
  \right )
  =
  -
  \sum _{ b _2,\,b _3,\,b _4,\,b _5 }
  \bar u _{ b _2 b _1 } u _{ b _3 b _1 } z _{ b _2 b _4 }
  s^{ b _4j } s ^{ i b _5 } \bar z _{ b _3 b _5 },
\end{eqnarray*} 
which completes the proof of the last assertion of this lemma. 
\hfill $\Box$

\begin{lem}\label{lem:a-5}
Let 
$ 
  \mathbf{\Phi} ( \mathbf{F} )
  =
  {\rm Diag\,}
  ( 
    \varphi _1 ( \mathbf{F} )
    ,\,
    \varphi _2 ( \mathbf{F} )
    ,\,
    \ldots
    ,\,
     \varphi _m ( \mathbf{F} )
   )
$, 
where 
$    \varphi _i ( \mathbf{F} )$'s $( i = 1,\,2,\,\ldots,\,m )$ are 
differentiable functions from $ \BB{R} ^m _{ > } \to \BB{R} _+ $. 
Then we have 
\begin{eqnarray*}
  {\rm Tr\,}\{{\rm Re\,}
  (
    \nabla _Z ^\prime
      \mathbf{U} \mathbf{\Phi} ( \mathbf{F} ) \mathbf{U} ^*
    \mathbf{Z}
  )
  \}
  &=&
  \sum _k 
  \biggl \{
    f _k \varphi _{ kk } ( \mathbf{F} )
    +
    (p-m+1) \varphi _k ( \mathbf{F} )
    +
    2 \sum _{ b > k }
    \frac{
      f _k \varphi _k ( \mathbf{F} )
      -
       f _{ b  } \varphi _{ b } ( \mathbf{F} )
    }{
      f _k - f _{ b }
    }
  \biggr \},
  \\
  {\rm Tr\,}
  ( 
    \mathbf{D} _S 
    \mathbf{Z} ^* \mathbf{U} \mathbf{\Phi} ( \mathbf{F} ) 
    \mathbf{U}^* \mathbf{Z}
  )
  &=&
  -
  \sum _k
  \biggl \{
   f _k ^2 \varphi _{ kk } ( \mathbf{F} )
    -2 ( m - 1 ) f _k \varphi _{k} ( \mathbf{F} ) 
    +
    2 \sum _{ b  > k }
    \frac{
      f _k ^2 \varphi _k ( \mathbf{F} )
      -
       f _b ^2 \varphi _b ( \mathbf{F} )
    }{
      f _k - f _b
    }
  \biggr \},
\end{eqnarray*}
where 
$
  \varphi _{ kk } ( \mathbf{F} ) 
  =
  (\partial \varphi _k / \partial f _k)( \mathbf{F} )
  ,\,
  k = 1,\,2,\,\ldots,\,m. 
$
\end{lem}
\par\noindent
{\bf Proof}. 
Use notation  $ \mathbf{\Phi}$, $ \varphi _k $, and   $ \varphi _{kk} $ 
short for $ \mathbf{\Phi} ( \mathbf{F} ) $, 
$ \varphi _{ k } ( \mathbf{F} )  $, and $
  \varphi _{ kk } ( \mathbf{F} ) $, respectively.  
To prove the first equation of this lemma, 
we first note that 
\[
  {\rm Tr\,}\{{\rm Re\,}(
  (
    \nabla _Z ^\prime
      \mathbf{U} \mathbf{\Phi} ( \mathbf{F} ) \mathbf{U} ^*
    \mathbf{Z}
  )
  \}
  =
  p {\rm Tr\,} (\mathbf{U} \mathbf{\Phi} \mathbf{U} ^*)
  +
  \frac{ 1 }{ 2 }
  \biggl ( 
    {\rm Tr\,} ( \mathbf{Z}  \nabla _Z ^\prime \mathbf{U} \mathbf{\Phi} \mathbf{U} ^* )
   +
      {\rm Tr\,} ( \bar{\mathbf{Z}}   \nabla _{ Z} ^*  \bar{\mathbf{U}} \mathbf{\Phi} \mathbf{U} ^\prime )
   \biggr ). 
\]
Now we use the first three equations in 
Lemma~$\ref{lem:a-4}$ and proceed in a  way similar to  the proof of Lemma~\ref{lem:a-3} in order to evaluate the second term inside the expectation of  the right hand side as 
\begin{eqnarray*}
  &&
  {\rm Tr\,} \{
    \mathbf{Z} 
    \nabla _Z ^\prime
      \mathbf{U} 
      \mathbf{\Phi} 
      \mathbf{U} ^* 
  \}
  =
  \sum _{ i,\,j,\,k,\,l }
  z _{ ij }
  \frac{
    \partial ( u _{ kl } \varphi _l \bar u _{ i l } )
  }{
    \partial z _{ kj }
  }
  \\
  &&
  \qquad
  =
  \sum _{ i,\,j,\,k,\,l }
  \biggl \{
  z _{ ij }
  \varphi _l \bar u _{ i l }
  \frac{
    \partial  u _{ kl } 
  }{
    \partial z _{ kj }
  }
  +
  z _{ ij } \varphi  _l u _{ kl }
  \frac{
    \partial \bar u _{ i l } 
  }{
    \partial z _{ kj }
  }
  +
  z _{ ij }
  u _{ kl }  \bar u _{ i l }
  \frac{
    \partial  \varphi _l 
  }{
    \partial z _{ kj }
  }
  \biggr \}
   \\
  &&
  \qquad
  =
  \sum _{ i,\,j,\,k,\,l,\,b _3,\,b_5}
   \biggl \{
  z _{ ij }
  \varphi _l \bar u _{ i l }
  \sum _{b _1 \not = l}
  \frac{
    u _{ k  b_1 }
    \bar u _{ k b _1 }
    u _{ b _3  l}
    s ^{ j b _5 }
    \bar z _{ b _3 b _5 }
  }{
    f _l - f _{ b _1 }
  }
  +
  z _{ ij } \varphi  _l u _{ kl }
  \sum _{ b _1 \not = l}
  \frac{
    \bar u _{ i b _1  }
    \bar u _{ k l }
    u _{b _3 b _1 }
    s ^{ j b _5 }
    \bar z _{ b _3 b _5 }
  }{
    f _l - f _{ b _1 }
  }
    \\
  &&
  \qquad\qquad
  +
  z _{ ij }
  u _{ kl }  \bar u _{ i l }
  \bar u _{  k b _1 }
  u _{ b _3 b _1 }
  s ^{ j b _5 }
  \bar z _{ b _3 b _5 }
  \frac{
    \partial  \varphi _l 
  }{
    \partial f _{ b _1 }
  }  
  \biggr \}
  \\
    &&
  \qquad
  =
  \sum _{ i,\,b _3,\,l }
  \biggl \{
  ( \mathbf{Z} \mathbf{S} ^{-1}\mathbf{Z} ^* ) _{ i b _3 }
  u _{ b _3 l }
  \sum _{ b _1 \not = l } 
  \frac{ \varphi _l }{ f _l - f _{ b _1 }}
  \bar u _{ i l}
  +
  ( \mathbf{Z} \mathbf{S} ^{-1}\mathbf{Z} ^* ) _{ i b _3 }
  u _{ b _3 b _1 }
   \sum _{ b _1 \not = l }  \frac{ \varphi _l }{ f _l - f _{ b _1} }
  \bar u _{ i b _1}  
  \\
     &&
  \qquad\quad
  +
  ( \mathbf{Z} \mathbf{S} ^{-1}\mathbf{Z} ^* ) _{ i b _3 }
  u _{ b _3 l }
  \frac{ \partial \varphi _l }{ \partial f _l }
  \bar u _{ i l}
  \biggr \}
  \\   
    &&
  \qquad
  =
  {\rm Tr\,} (
    \mathbf{Z} \mathbf{S} ^{-1}\mathbf{Z} ^*
    \mathbf{U} \widetilde{\mathbf{\Phi}} 
    \mathbf{U} ^*
  ),
\end{eqnarray*}
where 
$  
  \widetilde{\mathbf{\Phi}}
  =
  {\rm Diag\,}
  (
    \widetilde{\mathbf{\phi}} _1 
    ,\,
     \widetilde{\mathbf{\phi}} _2 
     ,\,
     \ldots
     ,\,
      \widetilde{\mathbf{\phi}} _m 
  )
$ 
with 
$
   \widetilde{\mathbf{\phi}} _i 
   =
   f _i
   (\partial \varphi _i/\partial f _i )
   +
   \sum _{ b \not = i }
   f _i 
   ( \varphi _i  - \varphi _b ) / ( f _i - f _b )
   ,\
   (i = 1,\,2,\,\ldots,\,m )
$. 
Similarly we use the first three equations in 
Lemma~$\ref{lem:a-4}$ to evaluate the third term inside the expectation of  the right hand side as 
\begin{eqnarray*}
  &&
  {\rm Tr\,} \{
    \bar{\mathbf{Z}}  
    \nabla _Z ^*
      \bar{\mathbf{U}} 
      \mathbf{\Phi} 
      \mathbf{U} ^\prime 
  \}
  =
  \sum _{ i,\,j,\,k,\,l }
  \bar z _{ ij }
  \frac{
    \partial ( \bar u _{ kl } \varphi _l  u _{ i l } )
  }{
    \partial \bar z _{ kj }
  }
  \\
  &&
  \qquad
  =
  \sum _{ i,\,j,\,k,\,l }
  \biggl \{
  \bar z _{ ij }
  \varphi _l  u _{ i l }
  \frac{
    \partial  \bar u _{ kl } 
  }{
    \partial \bar z _{ kj }
  }
  +
  \bar z _{ ij } \varphi  _l \bar u _{ kl }
  \frac{
    \partial  u _{ i l } 
  }{
    \partial \bar z _{ kj }
  }
  +
  \bar z _{ ij }
  \bar u _{ kl }   u _{ i l }
  \frac{
    \partial  \varphi _l 
  }{
    \partial \bar z _{ kj }
  }
  \biggr \}
  \\
  &&
  \qquad
  =
  \sum _{ i,\,j,\,k,\,l }
  \biggl \{
  \bar z _{ ij }
  \varphi _l  u _{ i l }
  \overline{
    \frac{
      \partial  u _{ kl } 
    }{
      \partial  z _{ kj }
    }
  }
  +
  \bar z _{ ij } \varphi  _l \bar u _{ kl }
  \overline{
    \frac{
      \partial  \bar u _{ i l } 
    }{
      \partial z _{ kj }
    }
  }
  +
  \bar z _{ ij }
  \bar u _{ kl }   u _{ i l }
  \overline{
    \frac{
      \partial  \varphi _l 
    }{
      \partial z _{ kj }
    }
  }
  \biggr \}
   \\
  &&
  \qquad
  =
  \sum _{ i,\,j,\,k,\,l,\, b _3,\,b_5}
   \biggl \{
  \bar z _{ ij }
  \varphi _l u _{ i l }
  \sum _{b _1 \not = l}
  \frac{
    \bar u _{ k  b_1 }
    u _{ k b _1 }
    \bar u _{ b _3  l}
    \bar s ^{ j b _5 }
    z _{ b _3 b _5 }
  }{
    f _l - f _{ b _1 }
  }
  +
  \bar z _{ ij } \varphi  _l \bar u _{ kl }
  \sum _{ b _1 \not = l}
  \frac{
    u _{ i b _1  }
    u _{ k l }
    \bar u _{b _3 b _1 }
    \bar s ^{ j b _5 }
    z _{ b _3 b _5 }
  }{
    f _l - f _{ b _1 }
  }
    \\
  &&
  \qquad\qquad
  +
  \bar z _{ ij }
  \bar u _{ kl }  u _{ i l }
  u _{ k b _1 }
  \bar u _{ b _3 b _1 }
  \bar s ^{ j b _5 }
  z _{ b _3 b _5 }
  \frac{
    \partial  \varphi _l 
  }{
    \partial f _{ b _1 }
  }  
  \biggr \}
  \\
    &&
  \qquad
  =
  \sum _{ i,\,b _3,\,l }
  \biggl \{
  ( \mathbf{Z} \mathbf{S} ^{-1}\mathbf{Z} ^* ) _{ b _3 i }
  u _{ il }
  \sum _{ b _1 \not = l }
  \frac{ \varphi _l }{ f _l - f _{ b _1 } }
  \bar u _{ b _3 l }
  +
  ( \mathbf{Z} \mathbf{S} ^{-1}\mathbf{Z} ^* ) _{ b _3 i }
  u _{ i b _1 }
   \sum _{ b _1 \not = l }
  \frac{ \varphi _l }{ f _l - f _{ b _1 } }
  \bar u _{ b _3 b _1}  
  \\
     &&
  \qquad\quad
  +
  ( \mathbf{Z} \mathbf{S} ^{-1}\mathbf{Z} ^* ) _{ b _3 i }
  u _{ i l }
    \frac{ \partial \varphi _l }{ \partial f _l }
  \bar u _{ b _3 l}
  \biggr \}
  \\   
    &&
  \qquad
  =
  {\rm Tr\,} (
    \mathbf{Z} \mathbf{S} ^{-1}\mathbf{Z} ^*
    \mathbf{U} \widetilde{\mathbf{\Phi}} 
    \mathbf{U} ^*
  ),
\end{eqnarray*}
Putting $ \mathbf{Z} \mathbf{S} ^{-1}\mathbf{Z} ^* = \mathbf{U}  \mathbf{F} \mathbf{U} ^* $ into the right hand side of the above two equations, we have 
\begin{eqnarray*}
    {\rm Tr\,}\{{\rm Re\,}
  (
    \nabla _Z ^\prime
      \mathbf{U} \mathbf{\Phi} ( \mathbf{F} ) \mathbf{U} ^*
    \mathbf{Z}
  )
  \}
  &=&
  \sum _k 
  \biggl \{
    p \varphi _k
    +
    f _k \varphi _{ kk }
    +
    \sum _{ b \not = k }
    \frac{
      f _k ( \varphi _k - \varphi _b )
    }{
      f _k - f _b
    }
  \biggr \}
  \\
  &=&
  \sum _k 
  \biggl \{
    p \varphi _k
    +
    f _k \varphi _{ kk }
    +
    \sum _{ b \not = k }
    \frac{
      f _k \varphi _k - f _b \varphi _b + (f _b-f _k ) \varphi _b 
    }{
      f _k - f _b
    }
  \biggr \} 
   \\
  &=&
  \sum _k 
  \biggl \{
    f _k \varphi _{ kk }
    +
    (p - m + 1)
    \varphi _k 
    +
    \sum _{ b \not = k }
    \frac{
      f _k \varphi _k - f _b \varphi _b  
    }{
      f _k - f _b
    }
  \biggr \}, 
\end{eqnarray*}
which completes the first equation of this lemma. 
\par
To prove the second equation of this lemma, we first note that
\begin{eqnarray}
 {\rm Tr\,}
  (
    \mathbf{D} _S \mathbf{Z} ^*
      \mathbf{U} \mathbf{\Phi}  \mathbf{U} ^*
    \mathbf{Z} 
  )
  &=&
  \sum _{ i,\,j,\,k _1,\,k _2,\,l }
  \frac{ 
  \partial ( 
       \bar z _{ k _1 j }
      u _{ k _1 l } 
      \varphi _l \bar u _{ k _2 l }
      z _{ k _2 i })
  }{
    \partial s _{ ij }
  }
  \nonumber
  \\
  &=&
  \sum _{ i,\,j,\,k _1,\,k _2,\,l }
  \bar z _{ k _1 j } z _{ k _2 i }
  \biggl \{
    \varphi _l \bar u_{ k _2 l }
    \frac{
      \partial u _{ k _1 l }
    }{
      \partial s _{ ij }
    }
    +
    \varphi _l u_{ k _1 l }
    \frac{
      \partial \bar u _{ k _2 l }
    }{
      \partial s _{ ij }
    }
    +
    u_{ k _1 l } \bar u _{ k _2 l }
    \frac{
      \partial \varphi _l
    }{
      \partial s _{ ij }
    }
  \biggr \}. 
 \label{eq:a-15}
\end{eqnarray}
But, from the forth equation of Lemma~\ref{lem:a-4}, we have 
\begin{eqnarray}
  \sum _{ i,\,j,\,k _1,\,k _2,\,l }
  \bar z _{ k _1 j } z _{ k _2 i }
    \varphi _l \bar u_{ k _2 l }
    \frac{
      \partial u _{ k _1 l }
    }{
      \partial s _{ ij }
    }
  &=&
  -\sum _{ b _1 \not = l }
  \frac{
    \varphi _l
    ( \mathbf{U} ^* \mathbf{Z} \mathbf{S} ^{-1} \mathbf{Z} ^* \mathbf{U} ) 
    _{ b _1 b _1 }
     ( \mathbf{U} ^* \mathbf{Z} \mathbf{S} ^{-1} \mathbf{Z} ^* \mathbf{U} ) 
    _{ ll }
  }{
    f _l - f _{ b _1 }
  }
  \nonumber
  \\
  &=&
  -\sum _{ b _1 \not = l }
  \frac{
    f _{ b _1 } f _l \varphi _l
  }{
    f _l - f _{ b _1 }
  },
  \nonumber
\end{eqnarray}
where we denote by 
$ ( \mathbf{U} ^* \mathbf{Z} \mathbf{S} ^{-1} \mathbf{Z} ^* \mathbf{U} ) 
    _{ ll }$ the $ l$-th diagonal element of a matrix 
    $ \mathbf{U} ^* \mathbf{Z} \mathbf{S} ^{-1} \mathbf{Z} ^* \mathbf{U} $ 
    for $ l = 1,\,2,\,\ldots,\,m$. 
Similarly, from the last two equations of Lemma~\ref{lem:a-4}, we have 
\begin{eqnarray}
  \sum _{ i,\,j,\,k _1,\,k _2,\,l }
  \bar z _{ k _1 j } z _{ k _2 i }
    \varphi _l u_{ k _1 l }
    \frac{
      \partial \bar u _{ k _2 l }
    }{
      \partial s _{ ij }
    }
  &=&
  -\sum _{ b _1 \not = l }
  \frac{
    f _{ b _1 } f _l \varphi _l
  }{
    f _l - f _{ b _1 }
  }
   \nonumber
\end{eqnarray}
and 
\begin{eqnarray}
  \sum _{ i,\,j,\,k _1,\,k _2,\,l }
  \bar z _{ k _1 j } z _{ k _2 i }
    u_{ k _1 l } \bar u _{ k _2 l }
    \frac{
      \partial \varphi _l
    }{
      \partial s _{ ij }
  }  
  &=&
  \sum _{ i,\,j,\,k _1,\,k _2,\,l,\,b _1 }
  \bar z _{ k _1 j }  z _{ k _2 i }
    u_{ k _1 l } \bar u _{ k _2 l }
    \frac{
      \partial \varphi _l
    }{
      \partial f _{ b _1 }
  }  
    \frac{
      \partial f _{ b _1 }
    }{
      \partial s _{ ij }
  }  
  \nonumber
  \\
  &=&
  -\sum _{ b _1,\,l }
    ( \mathbf{U} ^* \mathbf{Z} \mathbf{S} ^{-1} \mathbf{Z} ^* \mathbf{U} ) 
    _{ b _1 l }
     ( \mathbf{U} ^* \mathbf{Z} \mathbf{S} ^{-1} \mathbf{Z} ^* \mathbf{U} ) 
    _{ l b _1 }
    \frac{
      \partial \varphi _l
    }{
      \partial f _{ b _1 }
    }  
  \nonumber
  \\
  &=&
  -\sum _{ l }
    f ^2 _l 
    \frac{
      \partial \varphi _l
    }{
      \partial f _{ l }
    }.     
   \nonumber
\end{eqnarray}
Putting these three above equations 
into $(\ref{eq:a-15})$ we conclude that
\begin{eqnarray*}
    {\rm Tr\,}
  (
    \mathbf{D} _S \mathbf{Z}
      \mathbf{U} \mathbf{\Phi}  \mathbf{U} ^*
    \mathbf{Z} ^*
  )
  &=&
  - 
  \sum _l \left \{
    2 \sum _{ b  \not = l }
    \frac{
      f_{ b } f _l \varphi _l
    }{
      f _l - f _b
    }
    + 
    f _l ^2 
    \frac{ 
      \partial \varphi _l
    }{
      \partial f _l
    }
  \right \}
  \\
  &=&
  - 
  \sum _l \left \{
    2 \sum _{ b  \not = l }
    \frac{
      \{f _l ( f _b - f _l ) + f _l ^2 \} \varphi _l
    }{
      f _l - f _b
    }
    +
    f _l ^2 
    \frac{ 
      \partial \varphi _l
    }{
      \partial f _l
    }
  \right \}  
    \\
  &=&
  \sum _l
  \left \{
    2 ( m - 1 ) f _l \varphi _l
    - 2 \sum _{ b > l }
    \frac{
      f _l ^2  \varphi _l - f _b ^2 \varphi _b
    }{
      f _l - f _b
    }
    - 
    f _l ^2 
    \frac{ 
      \partial \varphi _l
    }{
      \partial f _l
    }
  \right \}
  ,  
\end{eqnarray*}
which completes the proof of the second assertion of this lemma.
\hfill $\Box$

\vskip 24pt
\begin{flushright}
\begin{tabular}{ll}
{\sc Faculty of Science}, 
{\sc Japan Women's University}
\\
2-8-1 Mejirodai, Tokyo 112-8681, Japan
\\
email: konno@fc.jwu.ac.jp
\end{tabular}
\end{flushright}
\end{document}